\documentclass[16pt]{article}
\usepackage{amsmath,amsfonts,amssymb}
\usepackage{verbatim}
\usepackage{latexsym}

\newtheorem{theorem}{Theorem}

\newtheorem{corollary}[theorem]{Corollary}

\newtheorem{definition}[theorem]{Definition}

\newtheorem{lemma}[theorem]{Lemma}

\newtheorem{remark}[theorem]{Remark}

\newenvironment{proof}[1][Proof]{\textbf{#1.} }{\ \rule{0.5em}{0.5em}}

\input{tcilatex}

\begin{document}

\title{Mixed Needlets}
\author{Daryl Geller and Domenico Marinucci \\
Department of Mathematics, Stony Brook University and \\
Department of Mathematics, University of Rome Tor Vergata}
\maketitle

\begin{abstract}
The construction of needlet-type wavelets on sections of the spin line
bundles over the sphere has been recently addressed in Geller and Marinucci
(2008), and Geller et al. (2008,2009). Here we focus on an alternative
proposal for needlets on this spin line bundle, in which needlet
coefficients arise from the usual, rather than the spin, spherical
harmonics, as in the previous constructions. We label this system \emph{%
mixed needlets} and investigate in full their properties, including
localization, the exact tight frame characterization, reconstruction
formula, decomposition of functional spaces, and asymptotic uncorrelation in
the stochastic case. We outline astrophysical applications.

\begin{itemize}
\item Keywords and phrases: Spherical Harmonics, Line Bundles, Spin
Needlets, Mixed Needlets, Besov Spaces, Cosmic Microwave Background
Radiation, Polarization, Weak Gravitational Lensing

\item AMS Classification: 42C40, 60G60, 33C55, 62M15, 83F05, 58J05
\end{itemize}
\end{abstract}

\section{Introduction and motivations}

A growing literature has been recently concerned with (random or
deterministic) spin functions, i.e. sections of spin fiber bundles over the
sphere (we defer a rigorous definition and more formal discussion to the
next section). Actually the interest of the physical literature on such
objects goes back for several decades, the seminal contributions going back
to \cite{newman} and \cite{goldb}. In these papers, spin spherical harmonics
were introduced in the language of physicists and used for the analysis of
gravitational radiation. Much more recently, spin functions have found a
crucial role in the analysis of cosmological observations, in particular in
connection with so-called Cosmic Microwave Background polarization data.
Polarization is a peculiar imprint characterizing the electromagnetic
radiation which was emitted at the age of recombination, some 13.7 billion
years ago and in the immediate adjacency of the Big Bang; as such, it
delivers information on a number of extremely important topics in the
current landscape of physical and cosmological research, for instance on the
existence of primordial gravitation waves, on the reionization era and on
primordial non-Gaussianity. The literature on these issues is vast; we
mention \cite{ck05,dode2004,selzad,kks96} for an introduction, while a
massive amount of observations are currently being collected by satellite
experiments by NASA and ESA (\emph{WMAP} and \emph{Planck}, respectively).
Spin fiber bundles will definitely be of the greatest interest for other
areas of physical research in the next decade, for instance in the analysis
of gravitational weak lensing on the images of galaxies \cite{bridle}. We
expect random sections of spherical fiber bundles to enjoy a growing
relevance also outside the physical sciences, for instance in medical
imaging.

Despite such a rich environment from the physical sciences, the interest in
the mathematical literature on spin bundles on the sphere has grown only
very recently. In particular, some efforts have been entertained to extend
to sections of spin fiber bundles the construction of spherical wavelets of
a needlet type. Scalar needlets were introduced for the sphere by \cite{npw1}%
, \cite{npw2}; the general case of compact Riemannian manifolds has
been presented by \cite{gm1,gm2,gm3}. The analysis of the asymptotic
properties of scalar needlets in random circumstances has been
started by \cite{bkmpAoS} ,\cite{bkmpBer}, see also \cite{ejslan},
\cite{bkmp09}, \cite{mayeli}, \cite{spalan}, \cite{bkmpAoSb},
\cite{kpp} for further developments and
\cite{pbm06},\cite{mpbb08},\cite{pietrobon1},
\cite{fay08},\cite{dela08},\cite{rudjord1},\cite{rudjord2},\cite{ghosh},
among others, for applications to cosmological data. Spin needlets
were introduced by \cite{gelmar}, in that paper, localization and
uncorrelation properties are also addressed. The general case where
the wavelet system need not be compactly supported in harmonic space
is discussed by \cite{gm4}; stochastic properties and related
statistical procedures are investigated in \cite{glm}, while
applications to a CMB framework are provided by \cite{ghmkp}. In
\cite{bkmp09}, it is proved that spin needlets actually make up a
tight frame system, with the same cubature points as the scalar
case, and characterizations of Besov spaces on spin fiber bundles
are discussed. Further results on the stochastic foundations of spin
random fields are provided by \cite{mal} and \cite{leosa}.

The purpose of this paper is to consider an alternative construction
for wavelets on spin fiber bundles. In particular, we focus on the
case where the resulting needlet coefficients are (complex-valued)
scalars, rather than spin quantities as in the spin needlet
proposal. We label this system \emph{mixed needlets} and investigate
in full their properties, including localization, the exact tight
frame characterization, reconstruction formula, decomposition of
functional spaces, and asymptotic uncorrelation in the stochastic
case; we outline also astrophysical applications. Concerning the
latter, we stress in particular that mixed needlets allow for
possibilities that were ruled out for the pure spin construction,
such as the estimation of cross power spectra between scalar and
spin components.

The plan of the paper is as follows: in Section 2 we review some background
material on spin fiber bundles, and we discuss some equivalent definitions
which have been provided in the literature. In Section 3 we explain the
construction of spin needlets, while Section 4 is devoted to the
investigation of mixed spin needlets. Section 5 discusses characterizations
of Besov spaces, establishing the equivalence of spin and mixed spin
needlets in this regard and investigating on the properties of functional
spaces for underlying scalar functions. Section 6 is devoted to directions
for future research, with particular reference to statistical applications.

\section{Spin functions}

\subsection{Some Definitions}

We begin by summarizing some definitions and basic facts about spin
functions. For more details, the reader may consult our article \cite{gelmar}%
.

We let $\mathbf{N}$ be the north pole of the unit sphere $\mathbb{S}^{2}$,
namely $(0,0,1)$, and we let $\mathbf{S}$ be the south pole, $(0,0,-1)$. Let
\begin{equation*}
U_{I}=\mathbb{S}^{2}\setminus \{\mathbf{N},\mathbf{S}\}.
\end{equation*}%
If $R\in SO(3)$, we define
\begin{equation*}
U_{R}=RU_{I}.
\end{equation*}%
On $U_{I}$ we use standard spherical coordinates $(\vartheta ,\varphi )$ ($%
0<\vartheta <\pi $, $-\pi \leq \varphi <\pi $), and analogously, on any $%
U_{R}$ we use coordinates $(\vartheta _{R},\varphi _{R})$ obtained by
rotation of the coordinate system on $U_{I}$.

At each point $p$ of $U_R$ we let $\rho_R(p)$ be the unit tangent vector at $%
p$ which is tangent to the circle $\vartheta_R = $ constant, pointing in the
direction of increasing $\varphi_R$. (This is well-defined since $R\mathbf{N}%
, R\mathbf{S} \notin U_R$.) We think of $\rho_R(p)$ as being our ``reference
direction'' at $p$, relative to the chart $U_R$. \newline
If $p \in U_{R_1} \cap U_{R_2}$, we let

\begin{center}
$\psi_{p R_2 R_1}$ be the oriented angle from the reference direction $%
\rho_{R_1}(p)$ to $\rho_{R_2}(p)$.
\end{center}

(See \cite{gelmar} for a precise explanation of which is the oriented angle.
For example, at $\mathbf{S}$, the oriented angle from $\vec{i}$ ti $\vec{j}$
is $\pi/2$; at $\mathbf{N}$, the oriented angle from $\vec{j}$ to $\vec{i}$
is $\pi/2$.)

Since $R$ is conformal, the angle $\psi_{p R_2 R_1}$ would be clearly be the
same if we had instead chosen $\rho_R(p)$ to point in the direction of
increasing $\vartheta_R$, for instance. Thus $\psi_{p R I}$ measures ``the
angle by which the tangent plane at $p$ is rotated if the coordinates are
rotated by $R$''.

Now say $\Omega \subseteq \mathbb{S}^{2}$ is open. Let\newline
\ \newline
$C_{s}^{\infty }(\Omega )=\{F=(F_{R})_{R\in SO(3)}:$ all $F_{R}\in C^{\infty
}(U_{R}\cap \Omega )$, and for all $R_{1},R_{2}\in SO(3)$ and all $p\in
U_{R_{1}}\cap U_{R_{2}}\cap \Omega $,
\begin{equation*}
F_{R_{2}}(p)=e^{is\psi }F_{R_{1}}(p),
\end{equation*}%
where $\psi =\psi _{pR_{2}R_{1}}\}$. \newline
\ \newline
(Heuristically, a physicist would think of $F_{R}$ as $F_{I}$ ``looked at
after the coordinates have been rotated by $R$''; at $p$, it has been
multiplied by $e^{is\psi }$, which is how physicists think of spin
quantities behaving after rotation.) Equivalently, $C_{s}^{\infty }(\Omega )$
consists of the smooth sections over $\Omega $ of the line bundle with
transition functions $e^{is\psi _{pR_{2}R_{1}}}$ from $U_{R_{1}}$ to $%
U_{R_{2}}$.

Similarly we can define $L_{s}^{2}(\Omega )$ (the $F_{R}$ need to be in $%
L^{2}$). There is a well-defined inner product on $L_{s}^{2}(\Omega )$,
given by $\langle F,G\rangle =\langle F_{R},G_{R}\rangle $; clearly this
definition is independent of choice of $R$.\newline
There is a unitary action of $SO(3)$ on $L_{s}^{2}(\mathbb{S}^{2})$, given
by $F\rightarrow F^{R}$, which is determined by the equation
\begin{equation}
(F^{R})_{I}(p)=F_{R}(Rp).  \label{rotact}
\end{equation}%
We think of $F^{R}$ as a ``rotate'' of $F$. Thus we have two important
relations: if we ``rotate coordinates'', we have
\begin{equation*}
F_{R_{2}}(p)=e^{is\psi }F_{R_{1}}(p),
\end{equation*}%
while if we ``rotate $F$'', we have
\begin{equation*}
(F^{R})_{I}(p)=F_{R}(Rp).
\end{equation*}%
Physicists would say that spin $s$ quantities need to be multiplied by a
factor $e^{is\psi }$ when rotated.

\subsection{Twisted bundles}

Although we will not make use of it in this article, it may be useful to
recast the previous discussion in a slightly different form, as suggested by %
\cite{mal}. View $SO(2)$ as a closed subgroup of $SO(3),$ with elements $%
k\in SO(2).$ This an Abelian subgroup, and assuming $k$ is parametrized by
the Euler angle $\gamma \in \lbrack 0,2\pi ]$ the irreducible
representations of $SO(2)$ are well-known to be one-dimensional and given by
$W_{s}(k):\mathbb{C}\rightarrow \mathbb{C},$ $W_{s}(k)x\mathbb{=}\exp
(is\gamma )x\mathbf{,}$ where $s\in \mathbb{N}$. Let $g\in SO(3),$ and
consider the action%
\begin{equation*}
k:\left\{ SO(3)\times \mathbb{C}\right\} \rightarrow \left\{ \mathbb{S}^{2},%
\mathbb{C}\right\} \text{ , }k(g,x)=(gk,\exp (is\gamma )x)\text{ .}
\end{equation*}%
We denote by $\mathcal{E}_{s}$ the quotient space of orbits of the above
action; that is, two elements $(g_{1},x_{1})$ and $(g_{2},x_{2})$ belong to
the same equivalence class if there exist $k\in SO(2)$ such that $%
(g_{2},x_{2})=(g_{1}k,W_{s}(k)x_{1}).$ For $s=0,$ this is clearly isomorphic
to $\left\{ \mathbb{S}^{2},\mathbb{C}\right\} ,$ i.e. the space of
complex-valued functions on the sphere. For $s\neq 0,$ we obtain indeed the
same spin fiber bundle we defined before $\left\{ \mathcal{E}_{s},\pi ,%
\mathbb{S}^{2}\right\} $, by taking the projection%
\begin{equation*}
\pi :\mathcal{E}_{s}\rightarrow \mathbb{S}^{2}\text{ , }\pi (g,x)=gSO(2)%
\text{ , }
\end{equation*}%
where we denoted as usual $gSO(2)$ the equivalence class $\left\{ gk,\text{ }%
k\in SO(2)\right\} ;$ for $g\in SO(3),$ this is isomorphic to $%
SO(3)/SO(2)\simeq \mathbb{S}^{2}.$

\subsection{Spin Spherical Harmonics}

\label{spsh}

Next we need some facts about the spin spherical harmonics; again the reader
may consult \cite{gelmar} for further details.

Let $f\in L^{2}($\textbf{$\mathbb{S}^{2}$}$\mathbf{)}$ the space of
square-integrable functions on the sphere; it is a well-known fact that the
following spectral representation holds, in the $L^{2}$ sense:%
\begin{equation*}
f(x)=\sum_{lm}a_{lm}Y_{lm}(x)\text{ , }a_{lm}=\int_{\mathbf{\mathbb{S}^{2}}%
}f(x)\overline{Y}_{lm}(x)dx\text{ ,}
\end{equation*}%
or more formally%
\begin{equation*}
L^{2}(\mathbf{\mathbb{S}^{2})=\bigoplus\limits_{l=0}^{\infty }\mathcal{H}_{l}%
\text{ ,}}
\end{equation*}%
where $\left\{ \mathcal{H}_{l}\right\} $ are the linear spaces spanned by
the standard spherical harmonics $\left\{ Y_{lm}:-l\leq m\leq l\right\} ,$
which are certain eigenfunctions of the (positive) spherical Laplacian%
\begin{equation*}
\Delta _{\mathbf{\mathbb{S}^{2}}}Y_{lm}=l(l+1)Y_{lm}\text{ .}
\end{equation*}%
Explicit expressions for the $Y_{lm}$ may be found in \cite{gelmar}.

We next define the \emph{spin-raising} operator $\eth$ and the \emph{%
spin-lowering} operator $\overline{\eth}$.

$\eth ,\overline{\eth }$ are maps which take smooth spin $s$ functions to
smooth spin $s+1$ (resp. $s-1$) functions, and which \emph{commute with the
actions of the rotation group} (\ref{rotact}). On a smooth spin $s$ function
$F$, we have $(\eth F)_{R}=\eth _{sR}F_{R}$, $(\overline{\eth }F)_{R}=%
\overline{\eth }_{sR}F_{R}$, where

\begin{eqnarray}
\eth _{sR}F_{R}(\vartheta ,\varphi ) &=&-\left( \sin \vartheta _{R}\right)
^{s}\left[ \frac{\partial }{\partial \vartheta _{R}}+\frac{i}{\sin \vartheta
_{R}}\frac{\partial }{\partial \varphi _{R}}\right] \left( \sin \vartheta
_{R}\right) ^{-s}F_{R}(\vartheta _{R},\varphi _{R})\text{ },  \label{edth} \\
\overline{\eth }_{sR}F_{R}(\vartheta _{R},\varphi _{R}) &=&-\left( \sin
\vartheta _{R}\right) ^{-s}\left[ \frac{\partial }{\partial \vartheta _{R}}-%
\frac{i}{\sin \vartheta _{R}}\frac{\partial }{\partial \varphi _{R}}\right]
\left( \sin \vartheta _{R}\right) ^{s}F_{R}(\vartheta _{R},\varphi _{R})%
\text{ }.  \label{baredth}
\end{eqnarray}

The spin $s$ spherical harmonics, defined for $l \geq |s|$, are then given
by
\begin{eqnarray*}
Y_{lm,s} &=&\left\{ \frac{(l-s)!}{(l+s)!}\right\} ^{1/2}(\eth )^{s}Y_{lm}%
\text{ , for }s>0\text{ ,} \\
Y_{lm,s} &=&\left\{ \frac{(l+s)!}{(l-s)!}\right\} ^{1/2}(-\overline{\eth }%
)^{-s}Y_{lm}\text{ , for }s<0\text{ ,}
\end{eqnarray*}

so that, if $l\geq |s|$,
\begin{eqnarray}
\eth Y_{lm,s} &=&\left[ \left( l-s\right) \left( l+s+1\right) \right]
^{1/2}Y_{lm,s+1}\text{ },  \label{spinh} \\
\overline{\eth }Y_{lm,s} &=&-\left[ \left( l+s\right) \left( l-s+1\right) %
\right] ^{1/2}Y_{lm,s-1}\text{ },  \label{spinhb}
\end{eqnarray}%
see also \cite{wiaux06} and \cite{mwen}. The $Y_{lm,s}$, for $l\geq |s|$, $%
-l\leq m\leq l$, form an orthonormal basis for $L_{s}^{2}$. In addition, $%
Y_{lm,s}$ is an eigenfunction of the (positive) \emph{spin spherical
Laplacian}
\begin{equation}
\Delta _{s}=%
\begin{cases}
-\overline{\eth }\eth & \mbox{ if }s\geq 0, \\
-\eth \overline{\eth } & \mbox{ if }s<0,%
\end{cases}
\label{delsdf}
\end{equation}%
acting on smooth spin functions, with eigenvalue
\begin{equation}
e_{ls}=(l-|s|)(l+|s|+1).  \label{evdf}
\end{equation}

If $s=0$, then $\Delta _{s}$ is just the usual (positive) spherical
Laplacian. The formal adjoint of $\eth $ (mapping smooth spin $s$ functions
on $\mathbb{S}^{2}$ to smooth spin $s+1$ functions) is $-\overline{\eth }$,
so that $\Delta _{s}$ is formally self-adjoint on smooth spin $s$ functions.

For $l\geq |s|$, we let $\mathcal{H}_{ls}$ denote the linear span of the $%
Y_{lm,s}$ for $-l\leq m\leq l$. $\mathcal{H}_{ls}$ is the eigenspace of $%
\Delta _{s}$ for the eigenvalue $e_{ls}$, and the direct sum of the $%
\mathcal{H}_{ls}$ (for $l\geq |s|$, $-l\leq m\leq l$) is all of $L_{s}^{2}(%
\mathbb{S}^{2})$.

We make some elementary observations about the eigenvalues $e_{ls}$.\newline
\ \newline
If $l \geq |s|$,
\begin{equation}  \label{elsz}
e_{ls} = l(l+1)-|s|(|s|+1) \leq l(l+1) = e_{l0};
\end{equation}

and for any $l, l^{\prime}$ (always nonnegative),
\begin{equation}  \label{elsad}
\sqrt{e_{l0}} + \sqrt{e_{l^{\prime},0}} < \sqrt{e_{l+l^{\prime}+1,0}};
\end{equation}

and for any $l \geq 0$ and any $s$,
\begin{equation}  \label{elsup}
e_{l0} \leq e_{l+|s|,s}.
\end{equation}

Here (\ref{elsz}) and (\ref{elsup}) are trivial. To prove (\ref{elsad}),
write $e_{l+l^{\prime}+1,0} = e_{l0} + e_{l^{\prime},0} +
2(l+1)(l^{\prime}+1)$, then square both sides of (\ref{elsad}), to see that
the inequality is equivalent to
\begin{equation*}
\sqrt{ll^{\prime}(l+1)(l^{\prime}+1)} < (l+1)(l^{\prime}+1),
\end{equation*}
which is evident.

Note that, for a spin $s$ function $F_{s}$, we may speak unambiguously of
the number $|F_{s}(x)|$, for any $x\in \mathbb{S}^{2}$. We now prove the
following inequality, by imitating the familiar method of proof in the case $%
s=0$:

\begin{lemma}
\label{linf} Say $l \geq |s|$, and that $Y \in \mathcal{H}_{ls}$. Then
\begin{equation}  \label{linfway}
\|Y\|_{\infty} \leq \sqrt{\frac{2l+1}{4\pi}} \|Y\|_2.
\end{equation}
In particular, for any $m$,
\begin{equation}  \label{linfwaysp}
\|Y_{lm,s}\|_{\infty} \leq \sqrt{\frac{2l+1}{4\pi}}.
\end{equation}
\end{lemma}

\textbf{Proof} Let $s^+ = \max(s,0)$.  In section 5 of \cite{gelmar} we showed that $Z_{l,s}$,
defined by
\begin{equation}
Z_{l,s}=(-1)^{s^{+}}\left[ \frac{2l+1}{4\pi }\right] ^{1/2}Y_{l,-s,s},
\label{szondfway}
\end{equation}%
which we called the $s$-zonal harmonic for $\mathcal{H}_{ls}$, has quite
similar properties to the usual zonal harmonic (for $\mathcal{H}_{l0}$).
Those properties, and (\ref{rotact}) imply that the following argument is
valid, just as in the case $s=0$: Say $p\in \mathbb{S}^{2}$, and choose $%
R\in SO(3)$ with $R\mathbf{N}=p$. Then
\begin{equation*}
|Y(p)|=|Y^{R}(N)|=|\langle Y^{R},Z_{l,s}\rangle |\leq \Vert Y^{R}\Vert
_{2}\Vert Z_{l,s}\Vert _{2}=\Vert Y\Vert _{2}\langle Z_{l,s},Z_{l,s}\rangle
^{1/2}=\Vert Y\Vert _{2}|Z_{l,s}(\mathbf{N})|^{1/2}.
\end{equation*}%
But
\begin{equation*}
|Z_{l,s}(\mathbf{N})|=|\sum_{m=-l}^{l}Y_{lm,s}(\mathbf{N})\overline{Y_{lm,s}(%
\mathbf{N})}|=\frac{2l+1}{4\pi };
\end{equation*}%
in fact, for any $x\in \mathbb{S}^{2}$, and any $R^{\prime }\in SO(3)$, one
has
\begin{equation}
\sum_{m=-l}^{l}Y_{lm,sR^{\prime }}(x)\overline{Y_{lm,sR^{\prime }}(x)}=\frac{%
2l+1}{4\pi }.  \label{krrsum}
\end{equation}%
This completes the proof.\newline
\ \newline
For $L\geq |s|$, let
\begin{equation}
V_{L,s}=\bigoplus_{l=|s|}^{L}\mathcal{H}_{ls}\text{ }.  \label{ghio}
\end{equation}
We note the following Bernstein-type lemma, adapted from \cite{gelpes}.

\begin{lemma}
\label{bern} A smooth spin $s$ function $F$ is in $V_{L,s}$ if and only if
there exist $A > 0$ and $B < e_{L+1,s}$ such that for every nonnegative
integer $N$,
\begin{equation}  \label{bernway}
\|(\Delta_s)^N F\|_2 \leq AB^N.
\end{equation}
\end{lemma}

\textbf{Proof} If $F \in V_{L,s}$, we surely have (\ref{bernway}), with $A =
\|F\|_2$, $B = e_{L,s}$, by the orthogonality of the eigenspaces $\mathcal{H}%
_{ls}$.

For the converse, say we have (\ref{bernway}). Suppose $l \geq L+1$, and $Y
\in \mathcal{H}_{l,s}$; it suffices to show that $\langle F, Y \rangle = 0$.
But, since $\Delta_s$ is formally self-adjoint, for any $N$ we have
\begin{equation*}
|\langle F, Y \rangle| = e_{ls}^{-N}|\langle F, \Delta_s^N Y \rangle| =
e_{ls}^{-N}|\langle \Delta_s^N F, Y \rangle| \leq A(\frac{B}{e_{ls}})^N
\|Y\|_2.
\end{equation*}
Since $B < e_{L+1,s} \leq e_{ls}$, this yields $\langle F, Y \rangle = 0$
upon letting $N$ go to infinity. This completes the proof.\newline
\ \newline
From this we find the following important \emph{product property}. (The case
$r = -s$ was first proved in \cite{bkmp09} by developing the ideas of the
subsection which follows. Here instead we adapt arguments from \cite{gelpes}%
.)

\begin{lemma}
\label{prod}
\begin{equation*}
V_{K,r}V_{L,s} \subseteq V_{K + L + |r+s|, r+s}.
\end{equation*}
\end{lemma}

\textbf{Proof} First note that the product of a smooth spin $r$ function
with a smooth spin $s$ function is a smooth spin $r+s$ function.

Note next that if $F,G$ are respectively smooth spin $p$ and spin $q$
functions, then
\begin{equation}
\eth (FG)=(\eth F)G+F(\eth G),  \label{dprod}
\end{equation}%
similarly for $\overline{\eth }$ in place of $\eth $. (\ref{dprod}) follows
at once from (\ref{edth}), once we note that for any $p,R$, as differential
operators
\begin{equation*}
\eth _{pR}=\eth _{0R}+p\cot \vartheta _{R},
\end{equation*}%
and similarly for $q$ or $p+q$ in place of $p$.

To prove the lemma, it suffices to show that if $|r| \leq k \leq K$ and $|s|
\leq l \leq L$, and if $F = Y_{k\mu,r}$ and $G = Y_{lm,s}$ for some $\mu, m$%
, then $FG \in V_{k + l + |r+s|, r+s}$.

Iterating (\ref{dprod}) and the companion equation for $\overline{\eth}$, we
find that we can write

\begin{equation}  \label{delN}
(-\Delta_s)^N (FG) = \sum_{j = 0}^{2N} T_j,
\end{equation}
where \newline

\begin{center}
each $T_j$ is a sum of $\binom{2N}{j}$ terms of the form $(D F) (D^{\prime}
G)$, where\\[0pt]
\ \\[0pt]
$D$ (resp. $D^{\prime}$) is a $j$-fold (resp. $(2N-j)$-fold) product of $%
\eth $'s and $\overline{\eth}$'s, in some order.
\end{center}

Note that the constants appearing on the right sides of (\ref{edth}) and (%
\ref{baredth}) are equal to $\pm \sqrt{e_{l,s}}$ or $\pm \sqrt{e_{l,-s}}$.
Thus, by (\ref{elsz}), if $D$ is as above,

\begin{equation*}
DF = 0 \mbox{ or } bY_{k\mu,r^{\prime}} \mbox{ for some } r^{\prime},
\mbox{
where } |b| \leq e_{k0}^{j/2}.
\end{equation*}

Thus, notation as above, by Lemma \ref{linf} we have

\begin{equation*}
\|(D F) (D^{\prime} G)\|_2 \leq \|DF\|_2 \|D^{\prime}G\|_{\infty} \leq \sqrt{%
\frac{2l+1}{4\pi}}e_{k0}^{j/2}e_{l0}^{(2N-j)/2},
\end{equation*}

so that, by the binomial theorem,

\begin{eqnarray*}
\|\Delta_s^N (FG)\|_2 & \leq & \sqrt{\frac{2l+1}{4\pi}}\sum_{j = 0}^{2N}%
\binom{2N}{j}e_{k0}^{j/2}e_{l0}^{(2N-j)/2} \\
& = & \sqrt{\frac{2l+1}{4\pi}}B_1^{2N}
\end{eqnarray*}

where
\begin{equation*}
B_1 = \sqrt{e_{k0}} + \sqrt{e_{l0}} < \sqrt{e_{k+l+1,0}}
\end{equation*}
by (\ref{elsad}). Set $B = B_1^2$. By (\ref{elsup}), $B <
e_{k+l+|r+s|+1,r+s} $. By Lemma \ref{bern}, we see that $FG \in V_{k + l +
|r+s|, r+s}$, as desired.

\subsection{Connection with Wigner D matrices}

Next, we shall explain the connection of spin spherical harmonics with
Wigner D matrices. This connection provides an alternative point of view,
but it is not necessary for the rest of the article. For further details on
this connection, the reader may consult \cite{glm} and \cite{mal}.

It is well-known that the elements $D_{m0}^{l},$ $m=-l,...,l$ of
Wigner's $D$ matrices are proportional to the standard spherical
harmonics $Y_{lm}$. It turns out that this equivalence holds in much
greater generality, in fact one has that (compare \cite{mal} for a
discussion of phase conventions)
\begin{eqnarray}
Y_{lm,s}(\vartheta ,\varphi ) &=&(-1)^{m^{+}}\sqrt{\frac{2l+1}{4\pi }}%
\overline{D_{m,-s}^{l}}(\varphi ,\vartheta ,0 ) \label{sphwig}
\text{ .}
\end{eqnarray}%
Here in place of $Y_{lm,s},\varphi ,\vartheta $, we should have written $%
Y_{lm,sR},\varphi _{R},\vartheta _{R}$ throughout, but we drop the reference
to the choice of chart for ease of notation whenever this can be done
without the risk of confusion.

Many of the properties of spin spherical harmonics follow easily from their
proportionality to elements of Wigner's $D$ matrices. Indeed, for instance,
their orthonormality
\begin{equation*}
\int_{\mathbf{\mathbb{S}^{2}}}Y_{lm,s}(p)\overline{Y_{l^{\prime }m^{\prime },s}%
}(p)dp=\int_{0}^{2\pi }\int_{0}^{\pi }Y_{lm,s}(\vartheta ,\varphi )\overline{%
Y_{l^{\prime }m^{\prime }}}(\vartheta ,\varphi )\sin \vartheta
d\vartheta d\varphi =\delta _{l}^{l^{\prime }}\delta _{m}^{m^{\prime
}}\text{ ;}
\end{equation*}%
is immediate. 
Also, viewing spin-spherical harmonics as functions on the group $SO(3)$
(i.e. identifying $p=(\vartheta ,\varphi )$ as the corresponding
rotation by means of Euler angles), and using (\ref{sphwig}) and the
group addition properties we obtain
easily%
\begin{eqnarray*}
\sum_{m=-l}^{l}Y_{lm,s}\left( p\right) \overline{Y_{lm,s}\left( p^{\prime
}\right) } &=&\frac{2l+1}{4\pi }\sum_{m}\overline{D_{m,-s}^{l}}(\varphi ,\vartheta ,0)%
D_{m,-s}^{l}(\varphi ^{\prime },\vartheta ^{\prime },0) \\
&=&\frac{2l+1}{4\pi }\overline{D_{s,-s}^{l}}(\psi (p,p^{\prime
}))\text{ ,}
\end{eqnarray*}%
where $\psi (p,p^{\prime })$ denotes the composition of the two rotations
(explicit formulae can be found in \cite{vmk}). In the special case $%
p=p^{\prime }$, we recover (\ref{krrsum}).

\subsection{E and M modes}

For a smooth spin function $F$ on $\mathbb{S}^{2}$, we have the expansion
\begin{equation}
F=\sum_{l}\sum_{m}a_{lm,s}Y_{lm,s}\text{ .}  \label{specspin}
\end{equation}%
with rapid decay of the $a_{lm,s}$ in $l$. From (\ref{specspin}), a further,
extremely important characterization of spin functions was first introduced
by \cite{newman}, see also \cite{edth} and \cite{gelmar} for a more
mathematically oriented treatment. In particular, it can be shown that there
exist a scalar complex-valued function
\begin{equation}
g(\vartheta ,\varphi )=\Re \left\{ g\right\} +i\Im \left\{ g\right\} \text{
, }  \label{nofun}
\end{equation}%
such that,%
\begin{equation*}
F_{s}=F^{E}+iF^{M}
\end{equation*}%
\begin{equation}
=\sum_{lm}a_{lm;E}Y_{lm,s}+i\sum_{lm}a_{lm;M}Y_{lm,s}\text{ ,}
\label{charspin}
\end{equation}%
where
\begin{equation*}
F^{E}=(\eth )^{s}\Re \left\{ g\right\} \text{ , }F^{M}=(\eth )^{s}\Im
\left\{ g\right\} \text{ .}
\end{equation*}%
Note that $a_{lm,s}=a_{lm;E}+ia_{lm;M}$, where $a_{lm;E}=\overline{a}%
_{l,-m;E},$ $a_{lm;M}=\overline{a}_{l,-m;M}.$ It is also readily seen that%
\begin{eqnarray*}
a_{lm,s}+\overline{a_{l,-m,s}}
&=&a_{lm;E}+ia_{lm;M}+a_{lm;E}-ia_{lm;M}=2a_{lm;E}\text{ ,} \\
a_{lm,s}-\overline{a_{l,-m,s}}
&=&a_{lm;E}+ia_{lm;M}-a_{lm;E}+ia_{lm;M}=2ia_{lm;M}\text{ .}
\end{eqnarray*}%
In the cosmological literature, $\left\{ a_{lm;E}\right\} $ and $\left\{
a_{lm;M}\right\} $ are labelled the $E$ and $M$ modes (or the electric and
magnetic components) of CMB\ polarization.

\section{Spin needlets}

We now recall the construction of spin needlets, see \cite{gelmar}, \cite%
{ghmkp}, \cite{gm4}, \cite{glm} and \cite{bkmp09} for further details and
discussions. Fix a ``dilation parameter'' $B>0$; $B$ is often chosen to be $2
$, but it is sometimes useful to let it take other values. Let $\phi $ be a $%
C^{\infty }$ function on $\mathbb{R},$ symmetric and decreasing on $\mathbb{R%
}^{+},$ supported in $\left| \xi \right| \leq 1,$ such that $0\leq \phi (\xi
)\leq 1$ and $\phi (\xi )=1$ if $\left| \xi \right| \leq \frac{1}{B}.$ Let

\begin{equation}  \label{bi}
b^{2}(\xi )=\phi (\tfrac\xi B)-\phi (\xi )\geq 0
\end{equation}

Note that supp$b \subseteq [1/B,B]$, and that

\begin{equation}  \label{bsum}
\sum_{j}b^{2}({\textstyle\frac{\xi}{B^{j}}})=\lim_{j\rightarrow \infty }\phi(%
{\textstyle\frac{\xi}{B^{j}}})=1\text{ for all } \xi > 0 \text{ .}
\end{equation}

Of course the sum on the left side of (\ref{bsum}) is zero if $\xi = 0$.

Let $T$ be a positive self-adjoint operator on a Hilbert space $\mathcal{H}$%
, and let $P$ be the projection onto the null space of $T$. It is a special
case of Theorem 2.1(b) of \cite{gm1}, that we may use the spectral theorem
to replace $\xi$ by $T$ in (\ref{bsum}), obtaining that

\begin{equation}  \label{bTsum}
\sum_{j}b({\textstyle\frac{\sqrt{T}}{B^{j}}})^2 = I-P \text{ .}
\end{equation}
where the sum converges strongly.

We take $\mathcal{H} = L^2_s$, $T = \Delta_s$, $P\mathcal{H} = \mathcal{H}%
_{|s|,s}$. Thus, if $F=\sum_{l}\sum_{m}a_{lm,s}Y_{lm,s}\in \mathcal{H}$, then

\begin{equation}
b({\textstyle\frac{\sqrt{T}}{B^{j}}})F =\sum_{l}\sum_{m}b(\frac{\sqrt{e_{ls}}%
}{B^{j}})a_{lm,s}Y_{lm,s}.  \label{bfirst}
\end{equation}

From this, it is easy to check (\ref{bTsum}) directly. Note also that $P%
\mathcal{H}$ is finite-dimensional (in fact, $2|s|+1$ dimensional). Note
moreover that $b({\textstyle\frac{\sqrt{T}}{B^{j}}}) \equiv 0$ for $j$
sufficiently negative, specifically if $B^{2j} < e_{|s|+1,s}$, the smallest
positive eigenvalue of $\Delta_s$.

For $x \in U_R$, let
\begin{equation}  \label{lamr}
\Lambda _{j}(x,y,R)=\sum_{l}\sum_m b(\frac{\sqrt{e_{ls}}}{B^{j}})Y_{lm,sR}(x)%
\overline{Y}_{lm,s}(y)\text{ .}
\end{equation}

Then evidently, if $F(y)=\sum_{l}\sum_{m}a_{lm,s}Y_{lm,s}(y)\in \mathcal{H}$%
, we have

\begin{equation}
\lbrack b({\textstyle\frac{\sqrt{T}}{B^{j}}})F]_{R}(x)=\sum_{l}\sum_{m}b(%
\frac{\sqrt{e_{ls}}}{B^{j}})a_{lm,s}Y_{lm,sR}(x)=\int \Lambda
_{j}(x,y,R)F(y)dy.  \label{blam}
\end{equation}

Here the integral is over $\mathbb{S}^{2}$. It is important to note that, in
the notation of subsection \ref{spsh},
\begin{equation}
b({\textstyle\frac{\sqrt{T}}{B^{j}}})F\in V_{L_{js},s},  \label{bvls}
\end{equation}%
where $L_{js}$ ($:=L_j$ if $s$ is understood) is the largest integer with $%
e_{L_{j}s}\leq B^{2(j+1)}$. In particular

\begin{equation}  \label{ljasy}
L_j \sim B^j
\end{equation}
as $j \to \infty$.

Now, take $F\in (I-P)\mathcal{H}$. Apply both sides of (\ref{bTsum}) to $F$,
and take the inner product with $F$. We find
\begin{equation}
\Vert F\Vert _{L_{s}^{2}}^{2}=\sum_{j}\Vert b({\textstyle\frac{\sqrt{T}}{%
B^{j}}})F\Vert _{L_{s}^{2}}^{2}=\sum_{j}\int |b({\textstyle\frac{\sqrt{T}}{%
B^{j}}})F|^{2}(x)dx\text{ ,}  \label{writeout}
\end{equation}

while, as long as $x \in U_R$,

\begin{equation}
|b({\textstyle\frac{\sqrt{T}}{B^{j}}})F|^{2}(x)=|\int \Lambda
_{j}(x,y,R)F(y)dy|^{2}\text{ .}  \label{writeoutlam}
\end{equation}

By (\ref{bvls}) and Lemma \ref{prod}, $|b({\textstyle\frac{\sqrt{T}}{B^{j}}}%
)F|^{2}\in V_{2L_{j},0}$. That is, it is the restriction to the
sphere of an ordinary polynomial of degree at most $2L_{j}\sim
B^{j}$. Accordingly, as noted by \cite{bkmp09} it is possible to
follow the method used in \cite{npw1} in the case $s=0$: By familiar
results for polynomials on the sphere, then, there is a constant
$c>0$, such
that for each $j$, there is a $c/(L_j+1)$-net \footnote{%
See e.g. \cite{bkmpBer} for the definition of $\varepsilon $-net.} $\{\xi
_{jk}\}$ of points on the sphere, and cubature weights $\{\lambda _{jk}\}
\sim (L_j+1)^{-2}$, such that for every polynomial $q$ of degree at most $2L_{j}$%
,
\begin{equation}
\int q(x)dx=\sum_{k}\lambda _{jk}q(\xi _{jk})  \label{cub}
\end{equation}

Thus, for $F\in (I-P)\mathcal{H}$, and provided $\xi _{jk}\in R_{jk}$, we in
fact have

\begin{equation}
\Vert F\Vert _{L_{s}^{2}}^{2}=\sum_{j}\sum_{k}\lambda _{jk}|\int \Lambda
_{j}(\xi _{jk},y,R_{jk})F(y)dy|^{2}\text{ .}  \label{frmsp}
\end{equation}

In other words, for $F\in (I-P)L_{s}^{2}$,

\begin{equation}
\Vert F\Vert _{L_{s}^{2}}^{2}=\sum_{j}\sum_{k}|\langle F,\psi _{jk,s}\rangle
|^{2},  \label{ttfrm}
\end{equation}

where the \emph{spin needlets} $\psi _{jk,s}$ are defined by

\begin{equation}  \label{psjk}
\psi _{jk,s}(y) =\sqrt{\lambda_{jk}}\: \overline{\Lambda}_j(%
\xi_{jk},y,R_{jk}) = \sqrt{\lambda _{jk}}\sum_{l}b(\frac{\sqrt{e_{ls}}}{B^{j}%
})\sum_{m}Y_{lm,s}(y)\overline{Y}_{lm,sR_{jk}}(\xi _{jk})\text{ .}
\end{equation}

Since $e_{|s|,s} = 0$, each $\psi _{jk,s} \in (I-P)L^2_s$. Consequently
\emph{the $\{\psi _{jk,s}\}$ are a tight frame for $(I-P)L^2_s$}.

For $F$ as above, we also define its \emph{spin needlet
coefficients} by
\begin{equation}
\beta _{jk,s}=\langle F,\psi _{jk,s}\rangle =\sqrt{\lambda _{jk}}\sum_{l}b(%
\frac{\sqrt{e_{ls}}}{B^{j}})\sum_{m}a_{lm,s}Y_{lm,sR_{jk}}(\xi _{jk})\text{ .%
}  \label{bejk}
\end{equation}

By general frame theory, if $F\in (I-P)L_{s}^{2}$, we have the
reconstruction formula

\begin{equation}
F=\sum_{j}\sum_{k}\beta _{jk,s}\psi _{jk,s}\text{ .}  \label{repfrm}
\end{equation}

\textbf{Remarks} 1. The choice of $R_{jk}$ does not affect any of the terms
on the right side of (\ref{ttfrm}) or (\ref{repfrm}). For this reason, and
for simplicity we will sometimes omit the $R_{jk}$ subscript in the formulas
(\ref{psjk}) and (\ref{bejk}) for $\psi _{jk,s}$ and $\beta _{jk,s}$, when
this can be done without causing confusion.\newline
2. We are ignoring the finite-dimensional space $PL_{s}^{2}$. This is
acceptable, because in astrophysical applications, the interest is in high
frequencies.\newline
3. One can use more general $b$, than those we used in (\ref{bi}), to
construct spin wavelets on the sphere, as in \cite{gelmar}, \cite{gm4}. This
leads to nearly tight frames with other interesting properties. For
instance, one can arrange for the support of the frame elements at scale $%
B^{-j}$ to be contained in a geodesic ball of radius $CB^{-j}$ (for some
fixed $C$).\newline
4. We will use the following notation and observations in Section 5. Let us
set
\begin{equation}
Q_{j}=b({\textstyle\frac{\sqrt{T}}{B^{j}}})^{2}.  \label{qjzer}
\end{equation}%
Following the arguments of (\ref{writeoutlam}) -- (\ref{ttfrm}), but now
without summing over $j$, we have that for $F\in L_{s}^{2}$,
\begin{equation}
\langle Q_{j}F,F\rangle =\sum_{k}|\langle F,\psi _{jk,s}\rangle |^{2},
\label{qjone}
\end{equation}%
After polarizing this identity, we see that for $F\in L_{s}^{2}$,
\begin{equation}
Q_{j}F=\sum_{k}\langle F,\psi _{jk,s}\rangle \psi _{jk,s}.  \label{qjtwo}
\end{equation}%
In (\ref{psjk}), $b(\frac{\sqrt{e_{ls}}}{B^{j}})=0$ unless $\frac{\sqrt{%
e_{ls}}}{B^{j}}\in (1/B,B)$. Thus, for all $j$,
\begin{equation}
\psi _{jk,s}\in V_{L_{j}}\ominus V_{L_{j-2}}.  \label{qjthree}
\end{equation}%
Accordingly,
\begin{equation}
Q_{j}:L_{s}^{2}\rightarrow V_{L_{j}}\ominus V_{L_{j-2}},\mbox{ and }%
Q_{j}\equiv 0\mbox{ on }[V_{L_{j}}\ominus V_{L_{j-2}}]^{\perp }.
\label{qjfour}
\end{equation}%
For any integer $N$, let
\begin{equation}
P_{N}=\sum_{j=-\infty }^{N}Q_{j}.  \label{qjfive}
\end{equation}%
As we know, the sum in (\ref{qjfive}) is actually finite, and, by (\ref%
{bTsum}), $P_{N}\rightarrow I-P$ strongly as $N\rightarrow \infty $. By
this, (\ref{qjthree}), and (\ref{qjfour}), we have
\begin{equation}
\psi _{jk,s}=(Q_{j-1}+Q_{j}+Q_{j+1})\psi _{jk,s},  \label{qjsix}
\end{equation}%
for all $j,k$.

\section{Mixed needlets}

We now present a construction of a different tight frame for spin functions,
which we shall call \emph{mixed needlets}.

For now we work with ordinary $L^{2}$ functions. Let $r$ be a fixed integer.
Let $\mathcal{H}_{l}$ be the space of spherical harmonics on $\mathbb{S}^{2}$
of degree $l$, and let $\mathcal{H}^{r}=\bigoplus_{l=|r|}^{\infty }\mathcal{H%
}_{l}\subseteq L^{2}$. We consider (\ref{bTsum}) with $\mathcal{H}=\mathcal{H%
}^{r}$ and with $T$ replaced by $\tilde{T}=\Delta -|r|(|r|+1)$, so that $P%
\mathcal{H}^{r}=\mathcal{H}_{|r|}$. Since $e_{lr}=e_{l0}-|r|(|r|+1)$, we
have that $\tilde{T}Y_{lm}=e_{lr}Y_{lm}$ for $l\geq |r|$.

In this situation, (\ref{bTsum}) leads to evident modifications of (\ref%
{bfirst})-- (\ref{frmsp}). To modify the equations, we simply replace $T$ by
$\tilde{T}$, and take $s = 0$, the sole exceptions being that we write $%
e_{lr}, e_{L_{jr}r}$ instead of $e_{l0}, e_{L_{j0}0}$. Of course we may
disregard all the rotations $R$ etc. We find that $\{\eta_{jk}\}$ is a tight
frame for $(I-P)\mathcal{H}$, where

\begin{equation}  \label{etjk}
\eta _{jk}(y)=\sqrt{\lambda _{jk}}\sum_{l = |r|}^{\infty} b(\frac{\sqrt{%
e_{lr}}}{B^{j}})\sum_{m}Y_{lm}(y)\overline{Y}_{lm}(\xi _{jk}).
\end{equation}

The last step in constructing mixed needlets is to change the variable name $%
r$ to $s$ in (\ref{etjk}), then to note that $\mathcal{H}^s$ is unitarily
equivalent to $L^2_s$, by means of the unitary equivalence $U$, where $%
U(Y_{lm}) = Y_{lm,s}$ (for $l \geq |s|$). Thus, if $P$ is now, once again,
as in the previous section, so that $PL^2_s = \mathcal{H}_{|s|,s}$, then the
\emph{mixed spin needlets} (or \emph{mixed needlets} for short) $\{\psi_{jk,s%
\mathcal{M}}\}$ are a tight frame for $(I-P)L^2_s$, where

\begin{equation}  \label{mpsjk}
\psi_{jk,s\mathcal{M}}(y)=\sqrt{\lambda _{jk}}\sum_{l = |s|}^{\infty} b(%
\frac{\sqrt{e_{ls}}}{B^{j}})\sum_{m}Y_{lm,s}(y)\overline{Y}_{lm}(\xi _{jk}).
\end{equation}

For $F=\sum_{l}\sum_{m}a_{lm,s}Y_{lm,s}\in L_{s}^{2}$, we also define its
\emph{spin needlet coefficients} by
\begin{equation}
\beta _{jk,s\mathcal{M}}=\langle F,\psi _{jk,s\mathcal{M}}\rangle =\sqrt{%
\lambda _{jk}}\sum_{l}b(\frac{\sqrt{e_{ls}}}{B^{j}})\sum_{m}a_{lm,s}Y_{lm}(%
\xi _{jk})\text{ .}  \label{mbejk}
\end{equation}

By general frame theory, if $F\in (I-P)L_{s}^{2}$, we have the
reconstruction formula

\begin{equation}
F=\sum_{j}\sum_{k}\beta _{jk,s\mathcal{M}}\psi _{jk,s\mathcal{M}}.
\label{mrepfrm}
\end{equation}

In this ``mixed'' situation, we set $Q_{j\mathcal{M}} = Ub({\textstyle\frac{%
\sqrt{\tilde{T}}}{B^{j}}})^2U^{-1}$. Then ``mixed'' analogues of (\ref{qjone}%
) -- (\ref{qjsix}) hold; it is only necessary to replace $Q_j$ by
$Q_{j\mathcal{M}}$ and $\psi _{jk,s}$ by $\psi _{jk,s\mathcal{M}}$
in those equations.\newline

But in fact we also have:

\begin{lemma}
\label{qjm}
\begin{equation*}
Q_{j\mathcal{M}} = Q_j.
\end{equation*}
\end{lemma}

\textbf{Proof} It suffices to show that these bounded operators agree on the
orthonormal basis elements $Y_{lm,s}$. But
\begin{equation*}
Q_{j\mathcal{M}}Y_{lm,s} = Ub({\textstyle\frac{\sqrt{\tilde{T}}}{B^{j}}}%
)^2U^{-1}Y_{lm,s} =
Ub({\textstyle\frac{\sqrt{\tilde{T}}}{B^{j}}})^2Y_{lm} =
Ub(\frac{\sqrt{e_{ls}}}{B^{j}})^2Y_{lm} = b(\frac{\sqrt{e_{ls}}}{B^{j}}%
)^2Y_{lm,s} = Q_jY_{lm,s}
\end{equation*}
as desired.\newline

In brief, the construction of spin needlets proceeds by applying the methods
of (\ref{bTsum}) -- (\ref{frmsp}) to $b({\textstyle\frac{\sqrt{T}}{B^{j}}})^2
$, while the construction of mixed needlets proceeds by applying the same
methods to the \emph{unitarily equivalent} operator $b({\textstyle\frac{%
\sqrt{\tilde{T}}}{B^{j}}})^2$, then invoking the unitary equivalence. Of
course this unitary equivalence is only effective on $L^2$, so there is no
evident reason to think that there would be an effective theory in other
function spaces for mixed needlets.  However, the main point of this article
is that mixed needlets do have nice mathematical properties beyond the $L^2_s
$ theory, as well as useful astrophysical applications. In particular, as we
shall show in the next section, they satisfy the usual needlet near-diagonal
localization property, in the same sense as spin needlets were shown to in %
\cite{gelmar}.

To understand better the meaning of mixed needlets and its relationship with
the existing literature, we start from (\ref{charspin}), and introduce the
notation%
\begin{equation}
f_{E}(x):= U^{-1}F^E = \sum_{lm}a_{lm;E}Y_{lm}(x)\text{ , }\:\:\:f_{M}(x):=
U^{-1}F^M = \sum_{lm}a_{lm;M}Y_{lm}(x)\text{ .}  \label{scalfun}
\end{equation}%
Clearly $f_{E}$ and $f_{M}$ are well-defined scalar functions which are
uniquely identified from $F_{s};$ as recalled above, in the $s=2$ of
interest for the physical literature they are labelled the electric and
magnetic components of the spin field (in the physical literature, $f_{M}$
is rather written $f_{B},$ but we already devoted the letter $B$ for another
purpose). Of course, it is possible to implement a standard (scalar) needlet
construction on these spaces, enjoying the well-known properties of needlets
(and indeed the same argument could be considered for other spherical
wavelets). The interesting question to address is clearly what are the
properties of such a procedure when viewed as acting on the original spin
space $L_{s}.$

More precisely, a direct idea to implement a wavelet transform on a spin
random field would be as follows. Start by to evaluating the spin transforms%
\begin{equation*}
\int_{\mathbb{S}^{2}}F_{s}(x)\overline{Y}_{lm,s}dx=a_{lm,s}\text{ , }\int_{%
\mathbb{S}^{2}}F_{s}(x)\overline{Y}_{lm,s}dx=a_{lm,s}\text{ ,}
\end{equation*}%
where
\begin{equation*}
a_{lm,s}=\frac{1}{2}\left\{ a_{lm,s}+\overline{a}_{l,-m,s}\right\} +\frac{1}{%
2}\left\{ a_{lm,s}-\overline{a}_{l,-m,s}\right\} =a_{lm;E}+ia_{lm;M}\text{ .%
}
\end{equation*}%
Note however that it is \emph{not} true that ${Re}(a_{lm;s})=a_{lm;E},$ ${Im}%
(a_{lm;s})=a_{lm;M},$ indeed $a_{lm;E},a_{lm;M}$ are complex-valued, and we
have%
\begin{eqnarray*}
a_{lm;E} &=&\frac{1}{2}\left\{ a_{lm,s}+\overline{a}_{l,-m,s}\right\} =\frac{%
1}{2}\left\{ a_{lm;E}+ia_{lm;M}+\overline{a}_{l,-m,E}-i\overline{a}%
_{l,-m,M}\right\} \text{ , } \\
a_{lm;M} &=&-\frac{i}{2}\left\{ a_{lm,s}-\overline{a}_{l,-m,s}\right\} =-%
\frac{i}{2}\left\{ a_{lm;E}+ia_{lm;M}-\overline{a}_{l,-m,E}+ia_{lm;M}\right%
\} \text{ ,}
\end{eqnarray*}%
where we use the (\emph{involutive}) property $a_{lm;E}=\overline{a}%
_{l,-m;E},$ $a_{lm;M}=\overline{a}_{l,-m;M}.$ This property uniquely
identifies the spherical coefficients $a_{lm;E},a_{lm;M}.$ Note that $%
a_{lm,s}$ is involutive if and only if the $M$ component is identically
null, while if and only if the $E$ component vanishes $ia_{lm,s}$ is
involutive.

It is then readily seen that%
\begin{eqnarray*}
\beta _{jk,\mathcal{M}} &:&=\int_{\mathbb{S}^{2}}F_{s}(x)\overline{\psi }%
_{jk,s\mathcal{M}}(x)dx=\sqrt{\lambda _{jk}}\sum_{lm}b(\frac{\sqrt{e_{ls}}}{%
B^{j}})a_{lm,s}Y_{lm}(\xi _{jk})
\end{eqnarray*}
\begin{equation*}
=\sqrt{\lambda _{jk}}\sum_{lm}b(\frac{\sqrt{e_{ls}}}{B^{j}})\left\{
a_{lm;E}+ia_{lm;M}\right\} Y_{lm}(\xi _{jk}) \\
\end{equation*}
\begin{equation}\label{ven}
=\beta _{jk;E}+i\beta _{jk;M},
\end{equation}
where $\beta _{jk;E},\beta _{jk;M}$ are real, i.e. $\beta
_{jk;E}={Re}(\beta
_{jk,\mathcal{M}}),$ $\beta _{jk;M}={Im}(\beta _{jk,\mathcal{M}}),$ because $a_{lm;E},a_{lm;M}$ and $%
Y_{lm}(\xi _{jk})$ are involutive. It is immediate to verify that
$\beta _{jk;E},\beta _{jk,M}$ are exactly the coefficients we would
obtain by evaluating a very slightly modified standard (scalar)
needlet transform on the scalar functions $f_{E},f_{M}.$ (It is
slightly modified because one would be using the $\eta_{jk}$ of
(\ref{etjk}) with $r=s$ instead of the
usual needlets, for which one would use $b(l/B^j)$ instead of $b(\sqrt{e_{ls}%
}/B^{j})$ in (\ref{etjk}).)

\subsection{Localization}

We fix an integer $s$. Because of (\ref{psjk}), localization properties of
the spin needlets $\psi_{jk,s}$ may be derived from localization properties
of the kernel $\Lambda_j(x,y,R)$, defined in (\ref{lamr}). Often we will
only need information about its absolute value $|\Lambda_j(x,y)|$ (which we
write as shorthand for $|\Lambda_j(x,y,R)|$ for any $R$).

For a full understanding, we need to consider the kernel

\begin{equation}  \label{lamtr}
\Lambda (x,y,t,R,R^{\prime},g)=\sum_{l}\sum_m g(t\sqrt{e_{ls}}) Y_{lm,sR}(x)%
\overline{Y}_{lm,sR^{\prime}}(y)\text{ ,}
\end{equation}

for $t > 0$, $R, R^{\prime}\in SO(3)$, $g \in C_c^{\infty}(\mathbb{R})$. In
particular

\begin{equation}  \label{lamrel}
|\Lambda_j (x,y)| = |\Lambda(x,y,B^{-j},R,R^{\prime},b)|,
\end{equation}

for any $R,R^{\prime}$.

When $g$ is fixed and understood we will write $\Lambda (x,y,t,R,R^{\prime})$
for $\Lambda (x,y,t,R,R^{\prime},g)$. In the case $s=0$, $\Lambda$ does not
depend on $R,R^{\prime}$, and we simply write $\Lambda(x,y,t)$ for $%
\Lambda(x,y,t,R,R^{\prime})$.

In the case $s = 0$, localization properties of a variant of this kernel
(where $e_{l0}$ in (\ref{lamtr}) is replaced by $l^2$) were derived in \cite%
{npw1}, \cite{npw2}. For the kernel as it stands in (\ref{lamr}), as well as
analogous kernels on smooth compact Riemannian manifolds, localization
results (including results where $b$ need not have compact support away from
$0$) were proved in \cite{gm1}.

For general $s$, the localization properties of $\Lambda$ were proved in %
\cite{gelmar} and \cite{gm4}.

For the purposes of this article, the relevant localization results are:

\begin{lemma}
\label{manmo} Say $s=0$. \newline
(a) Suppose $g\in C_{c}^{\infty }(0,\infty )$. Then for every pair of $%
C^{\infty }$ differential operators $X$ $($in $x)$ and $Y$ $($in $y)$ on $%
\mathbb{S}^{2}$, and for every integer $\tau \geq 0$, there exists $C>0$ as
follows. Suppose $\deg X=j$ and $\deg Y=k$. Then
\begin{equation}
\left| XY\Lambda (x,y,t)\right| \leq \frac{Ct^{-2-j-k}}{\left\{ 1+\frac{%
d(x,y)}{t}\right\} ^{\tau }}\text{,}  \label{diagest}
\end{equation}%
for all $t>0$ and all $x,y\in \mathbb{S}^{2}$.\newline
(b) If $A>0$ is fixed, and $g\in C_{c}^{\infty }(\mathbb{R})$ is even, then
the conclusion of (a) remains true as long as $0<t<A$.
\end{lemma}

\textbf{Remarks} 1. In relation to the function $f$ of \cite{gm1}, our $%
g(\xi) = f(\xi^2)$.\newline
2. Part (b) was shown in \cite{gm1} for $0 < t < 1$ -- see the last
paragraph of section 4 of that article. For $t \in [1,A]$, the estimate (\ref%
{diagest}) is trivial, for then the right side is uniformly bounded below by
a positive constant, and the left side is uniformly bounded above by a
positive constant, as is apparent from (\ref{lamtr}).)\newline
3. Say $[p,q] \subseteq (0,\infty)$. Using the remarks preceding Theorem 6.1
in \cite{gelmar}, one sees that the constant $C$ appearing in (\ref{diagest}%
) may be taken to be uniform for $g$ ranging over any bounded subset of the
Fr\'echet space $C_c^{\infty}([p,q])$. Indeed, this follows easily from
Lemma \ref{manmo} and the closed graph theorem.

\begin{lemma}
\label{manmols} Say $s$ is a fixed integer. \newline
(a) \cite{gelmar} Suppose $g \in C_c^{\infty}(0,\infty)$. Say $%
R,R^{\prime}\in SO(3)$. Then for every pair of compact sets $\mathcal{F}_R
\subseteq U_R$ and $\mathcal{F}_{R^{\prime}} \subseteq U_{R^{\prime}}$,
every pair of $C^{\infty}$ differential operators $X$ $($in $x)$ on $U_R$
and $Y$ $($in $y)$ on $U_{R^{\prime}}$, and for every integer $\tau \geq 0$,
there exists $C > 0$ as follows. Suppose $\deg X = j$ and $\deg Y = k$. Then
\begin{equation}  \label{diagests}
\left|XY\Lambda(x,y,t,R,R^{\prime})\right| \leq \frac{C t^{-2-j-k}}{\left\{
1+\frac{d(x,y)}{t}\right\} ^{\tau }}\text{,}
\end{equation}
for all $t > 0$, all $x \in \mathcal{F}_R$ and all $y \in \mathcal{F}%
_{R^{\prime}}$.\newline
(b) \cite{gm4} If $A > 0$ is fixed, and $g \in C_c^{\infty}(\mathbb{R})$ is
even, then the conclusion of (a) remains true as long as $0 < t < A$.
\end{lemma}

\textbf{Remark} Note that if $X=Y=$ the identity operator, $\left| XY\Lambda
(x,y,t,R,R^{\prime })\right| $ is independent of $R,R^{\prime }$, and then (%
\ref{diagests}) holds for all $x,y\in \mathbb{S}^{2}$.\newline

We are now going to show that a similar localization result to Lemma \ref%
{manmols} (a) holds for the mixed needlets.

Consider then the ``mixed needlet kernel''
\begin{equation*}
\Lambda _{\mathcal{M}}(x,y,t,R)=\sum_{lm}b(t\sqrt{e_{ls}})Y_{lm,sR}(x)%
\overline{Y}_{lm}(y)\text{ ;}
\end{equation*}%
It should be observed that
\begin{equation}
\Lambda _{\mathcal{M}}(x,x,t)\equiv 0\text{ for all }x\in \mathbb{S}^{2},%
\text{ }t>0,  \label{diagzero}
\end{equation}%
because%
\begin{equation*}
\left| \sum_{lm}b(t\sqrt{e_{ls}})Y_{lm,sR}(x)\overline{Y}_{lm}(x)\right|
=\left| \sum_{l}b(t\sqrt{e_{ls}})\sum_{m}D_{msR}^{l}(x)\overline{D}%
_{m0}^{l}(x)\right| =0\text{ ,}
\end{equation*}%
by the unitarity properties of the Wigner's $D$ matrices. In contrast, in
the unmixed situation $|\Lambda (x,x,t)|=\frac{2l+1}{4\pi }|\sum_{l}b(t\sqrt{%
e_{ls}})|$, by (\ref{krrsum}). If $0\neq b\geq 0$, as is usually the case in
applications, this is nonzero. Qualitatively, then, $\Lambda _{\mathcal{M}}$
is different from $\Lambda $; (\ref{diagzero}) might even lead one to
suspect that $\Lambda _{\mathcal{M}}$ might not be well-localized. However,
the methods of our article \cite{gelmar} show that it is:

\begin{lemma}
\label{manmolm} Say $s$ is a fixed integer. \newline
(a) Suppose $b\in C_{c}^{\infty }(0,\infty )$. Say $R,R^{\prime }\in SO(3)$.
Then for every pair of compact sets $\mathcal{F}_{R}\subseteq U_{R}$ and $%
\mathcal{F}_{R^{\prime }}\subseteq U_{R^{\prime }}$, every pair of $%
C^{\infty }$ differential operators $X$ $($in $x)$ on $U_{R}$ and $Y$ $($in $%
y)$ on $\mathbb{S}^{2}$, and for every integer $\tau \geq 0$, there exists $%
C>0$ as follows. Suppose $\deg X=j$ and $\deg Y=k$. Then
\begin{equation}
\left| XY\Lambda _{\mathcal{M}}(x,y,t,R)\right| \leq \frac{C_{\tau
}t^{-2-j-k}}{\left\{ 1+\frac{d(x,y)}{t}\right\} ^{\tau }}\text{,}
\label{diagestm}
\end{equation}%
for all $t>0$, all $x\in \mathcal{F}_{R}$ and all $y\in \mathbb{S}^{2}$.
\end{lemma}

\begin{proof}
We shall modify the proof of Theorem 6.1 of \cite{gelmar}. To make that
easier, we set $f(\xi^2) = b(\xi)$, so that $f \in C_c^{\infty}(0,\infty)$;
say that in fact supp$f \subseteq [p,q]$, where $p > 0$. We assume $s > 0$;
similar arguments will apply when $s < 0$. Of course the case $s = 0$ is
handled by Lemma \ref{manmo}.

We have
\begin{equation*}
\Lambda _{\mathcal{M}}(x,y,t,R)=\sum_{l\geq s}f(t^{2}e_{ls})\mathcal{K}_{%
\mathcal{M}R}^{l}(x,y)\text{ ,}
\end{equation*}%
where%
\begin{equation*}
\mathcal{K}_{\mathcal{M}R}^{l}(x,y)=\sum_{m}Y_{lm,sR}(x)\overline{Y}_{lm}(y)%
\text{ .}
\end{equation*}%
Thus%
\begin{equation*}
\Lambda _{\mathcal{M}}(x,y,t,R)=\eth _{Rx}^{[s]}\sum_{l\geq s}f(t^{2}e_{ls})%
\sqrt{\frac{(l-s)!}{(l+s)!}}\mathcal{K}_{0}^{l}(x,y)\text{ ,}
\end{equation*}%
where
\begin{equation*}
\mathcal{K}_{0}^{l}(x,y)=\sum_{m}Y_{lm}(x)\overline{Y}_{lm}(y)\text{ .}
\end{equation*}%
and here
\begin{equation*}
\eth _{Rx}^{[s]}=\eth _{s-1,R}\circ \ldots \circ \eth _{0,R}
\end{equation*}%
in the $x$ variable.

As in the proof of Theorem 6.1 of \cite{gelmar}, we note that
\begin{equation}
\frac{(l+s)!}{(l-s)!}=\prod\limits_{k=1}^{s}\left[ e_{l0}-\gamma _{k}\right]
\text{ , }  \label{facbrk}
\end{equation}%
where $\gamma _{k}:=k(k-1)$. As in \cite{gelmar} we choose $T_{0},T_{1}>0$
with $\gamma _{s}T_{0}^{2}<p/2$, $e_{s+1,s}T_{1}^{2}>q$. We note that $%
\Lambda _{\mathcal{M}}(x,y,t,R)\equiv 0$ for $t\geq T_{1}$, and that (\ref%
{diagestm}) is trivial for $t$ in the compact interval $[T_{0},T_{1}]\subset
(0,\infty )$. (Indeed, there the right side of (\ref{diagestm}) is uniformly
bounded below by a positive constant, and the left side is uniformly bounded
above by a positive constant.) It is then enough to focus on $t\in
(0,T_{0}], $ We now define
\begin{equation*}
f_{t}(u):=\frac{f(u-s(s+1)t^{2})}{\sqrt{\prod\limits_{k=1}^{s}\left[
u-\gamma _{k}t^{2}\right] }},
\end{equation*}%
supported in the fixed compact interval $[p,q_{1}]:=[p,q+s(s+1)T_{0}^{2}]$;
we note that the denominator does not vanish in the interval. Then, using (%
\ref{elsz}) and (\ref{facbrk}), we write%
\begin{equation*}
\sum_{l\geq s}f(t^{2}e_{ls})\sqrt{\frac{(l-s)!}{(l+s)!}}\mathcal{K}%
_{0}^{l}(x,y)=t^{s}\sum_{l\geq s}f_{t}(t^{2}e_{l0})\mathcal{K}%
_{0}^{l}(x,y)=t^{s}\Lambda _{\lbrack t]}(x,y)\text{ ,}
\end{equation*}%
for%
\begin{equation*}
\Lambda _{\lbrack t]}(x,y)=\sum_{l\geq s}f_{t}(t^{2}e_{l0})\mathcal{K}%
_{0}^{l}(x,y)\text{ .}
\end{equation*}%
Now note that the functions $f_{t}$ for $t\in (0,T_{0}]$ form a bounded
subset of $C_{c}^{\infty }([p,q_{1}])$, and recall Remark 3, after Lemma \ref%
{manmo}. Choose smooth differential operators on $\mathbb{S}^{2}$, in $x$,
of degree $s$ and $j$, which agrees with $\eth _{Rx}^{[s]}$ and $X$
respectively, in a neighborhood of $\mathcal{F}_{R}$. As in the proof of
Theorem 6.1 of \cite{gelmar}, we find
\begin{equation*}
\left| XY\Lambda _{\mathcal{M}}(x,y,t,R)\right| =t^{s}\left| XY\eth
_{Rx}^{[s]}\Lambda _{\lbrack t]}(x,y)\right| \leq ct^{s}\frac{Ct^{-2-j-s-k}}{%
\left\{ 1+d(x,y)/t\right\} ^{\tau }}=\frac{Ct^{-2-j-k}}{\left\{
1+d(x,y)/t\right\} ^{\tau }}\text{ ,}
\end{equation*}%
as desired.
\end{proof}

\begin{remark}
For astrophysical applications, it is very common to observe spin random
fields only in a subset of the sphere, i.e. $\mathbb{S}^{2}\backslash G,$
say, where $G\subset \mathbb{S}^{2}$ is a region contaminated by foreground
emission, for instance the Milky Way radiation. Lemma \ref{manmolm} implies
that, for $\tau =1,2,...$
\begin{eqnarray*}
\left| \beta _{jk;E}-{Re}\left\{ \int_{\mathbb{S}^{2}\backslash G}F_{s}(x)%
\overline{\psi }_{jk,s\mathcal{M}}dx\right\} \right|  &=&\left| {Re}\left\{
\beta _{jk;E}-\int_{\mathbb{S}^{2}\backslash G}F_{s}(x)\overline{\psi }_{jk,s%
\mathcal{M}}dx\right\} \right|  \\
&\leq &\left| \int_{G}F_{s}(x)\overline{\psi }_{jk,s\mathcal{M}}dx\right|  \\
&\leq &\int_{G}\left| F_{s}(x)\right| \left| \psi _{jk,s\mathcal{M}}\right|
dx \\
&\leq &\left\{ \sup_{x\in G}\left| F_{s}(x)\right| \right\} \frac{C_{\tau
}\mu (G)}{\left\{ 1+B^{j}d(\xi _{jk},G)\right\} ^{\tau }}\text{ ,}
\end{eqnarray*}%
\begin{eqnarray*}
\left| \beta _{jk;M}-{Im}\left\{ \int_{\mathbb{S}^{2}\backslash G}F_{s}(x)%
\overline{\psi }_{jk,s\mathcal{M}}dx\right\} \right|  &=&\left| {Im}\left\{
i\beta _{jk;M}-\int_{\mathbb{S}^{2}\backslash G}F_{s}(x)\overline{\psi }%
_{jk,s\mathcal{M}}dx\right\} \right|  \\
&\leq &\left| \int_{G}F_{s}(x)\overline{\psi }_{jk,s\mathcal{M}}dx\right|  \\
&\leq &\left\{ \sup_{x\in G}\left| F_{s}(x)\right| \right\} \frac{C_{\tau
}\mu (G)}{\left\{ 1+B^{j}d(\xi _{jk},G)\right\} ^{\tau }}\text{ ,}
\end{eqnarray*}%
where $\mu (.)$ denotes Lebesgue measure on $\mathbb{S}^{2}.$ In other
words, for all the coefficients corresponding to locations $\xi _{jk}\in
\mathbb{S}^{2}\backslash G_{\varepsilon }$ (where $G_{\varepsilon }:=\xi \in
\mathbb{S}^{2}:d(\xi ,G)>\varepsilon )$ the mixed needlet coefficients are
asymptotically unaffected by the presence of unobserved regions. This is
clearly a property of the greatest importance for cosmological applications.
\end{remark}

We now derive the following corollaries from Lemmas \ref{manmols} and \ref%
{manmolm}, which will be essential for the characterizations of Besov spaces
in the next section.

\begin{corollary}
\label{loccor1} (a) In (a) or (b) of Lemma \ref{manmols}, we have that for
some $C > 0$,
\begin{equation*}
\int |\Lambda(x,y,t)| dx \leq C,\:\:\: \int |\Lambda(x,y,t)|dy \leq C,
\end{equation*}
where the integrals are over $\mathbb{S}^2$.\newline
(b) In Lemma \ref{manmolm},  we have that for some $C > 0$,
\begin{equation*}
\int |\Lambda_{\mathcal{M}}(x,y,t)| dx \leq C,\:\:\: \int |\Lambda_{\mathcal{%
M}}(x,y,t)|dy \leq C,
\end{equation*}
where the integrals are over $\mathbb{S}^2$.
\end{corollary}

\textbf{Proof} This follows at once from the fact that $\int [1 +
d(x,y)/t]^{-N}dx \leq C_Nt^n$ for any $N > n$ (see for example, the third
bulleted point after Proposition 3.1 of \cite{gm1}).

\begin{corollary}
\label{allgood} Say $1 \leq p \leq \infty$.\newline{}(a) For each
$j$, $Q_j, P_j: C^{\infty}_s \to L_s^p$.  If $p \geq 2$, the
restriction of $Q_j, P_j$ to $L_s^p \subseteq L_s^2$ is bounded on 
$L_s^p$.  If $p < 2$, $Q_j, P_j$ may be extended from $C^{\infty}_s$
to be bounded operators on $L_s^p$. 

Further, the operators $Q_j, P_j$ are
uniformly bounded on $L_s^p$ for $-\infty < j <
\infty$.\newline{}(b) $\left\| \psi _{jk,s}\right\|_1, \left\| \psi
_{jk,s\mathcal{M}}\right\|_1 \leq C2^{-j}$.
\end{corollary}

\textbf{Proof} For (a), recall from (\ref{qjzer}) that $Q_j = b({\textstyle%
\frac{\sqrt{T}}{B^{j}}})^2$ on $L^2_s$.%
Take $\Lambda$ as in (\ref{lamtr}) for $g = b^2$. Then for $x
\in U_R$, $y \in U_{R^{\prime}}$, we have
\begin{equation}  \label{qjform}
(Q_jF)_R(x) = \int \Lambda (x,y,t,R,R^{\prime}) F_{R^{\prime}}(y)dy
\end{equation}
so that
\begin{equation*}
|(Q_jF)(x)| \leq \int |\Lambda (x,y,t)| |F(y)| dy
\end{equation*}
for $F \in L^2_s$. Part (a) for $Q_j$ is apparent from this and from
Corollary \ref{loccor1}(a). (Note also that (\ref{qjform}) continues to hold
for all $F \in L^p_s$ (one uses a density argument if $p < 2$).) Similarly, part (a) for $P_j$
also follows from Corollary \ref{loccor1} (a) (where one references part (b)
of Lemma \ref{manmols}), because $P_j = \phi({\textstyle\frac{\sqrt{T}}{%
B^{j+1}}})$ on $L^2_s$ by (\ref{bi}). (Note that $\phi$ equals $1$ near $0$,
and so has an even extension to a $C_c^{\infty}$ function.) \newline
Finally, part (b) follows from Corollary \ref{loccor1} (b), once we observe
that
\begin{equation}  \label{psimlam}
|\psi _{jk,s\mathcal{M}}(x)| = \sqrt{\lambda_{jk}}|\Lambda_{\mathcal{M}%
}(x,\xi_{jk},B^{-j})|,\:\:\:\:\:\lambda_{jk} \sim B^{-2j}.
\end{equation}
This completes the proof.\newline

For notational simplicity, we take $B=2$ for the rest of this article; the
results would easily generalize to general $B$.

Using Lemma \ref{manmolm}, we obtain the following estimates on the $%
L_{s}^{p}$ norms of the $\psi _{jk,s\mathcal{M}}$.

\begin{lemma}
\label{lpnrm} For $1 \leq p \leq \infty$, we have
\begin{equation}  \label{lpnrmway}
\left\| \psi _{jk,s\mathcal{M}}\right\|_p \sim
2^{j(1-\frac{2}{p})}.
\end{equation}
\end{lemma}

\textbf{Proof} Let us call the estimate $\left\| \psi
_{jk,s\mathcal{M}}\right\|_p \leq C2^{j(1-\frac{2}{p})}$ the
majorization for this value of $p$, and call the reverse inequality
the minorization.

First we do the cases $p=1,2,\infty$. For $p = \infty$, the majorization
(by $C2^j$) follows at once from (\ref{psimlam}) and Lemma \ref{manmolm}. For $p = 1$, the
majorization (by $C2^{-j}$) is Corollary \ref{allgood} (b). For $p = 2$, the
majorization (by a constant) follows at once from the tight frame property.

For the remaining estimates, we adapt arguments from \cite{bkmp09} and \cite%
{bkmpBer}. Fix $c > 0$, $0 < \nu < 1$ such that $b^2 > c$ on the interval $%
[\nu,1]$. Using the orthonormality of the spin spherical harmonics, one
obtains the minorization for $p = 2$ from

\begin{equation*}
\left\| \psi _{jk,s\mathcal{M}}\right\| _2^{2} = \sum_{lm}\lambda_{jk}
\left| Y_{lm}(\xi _{jk})\right| ^{2}b^{2}(\frac{\sqrt{e_{ls}}}{2^{j}})\geq
c2^{-2j}\sum_{\{l:\: \nu^2 \leq e_{ls}/4^j \leq 1\}} \frac{2l+1}{4\pi} c
\:\:\geq\:\: c^{\prime }>0\text{ .}
\end{equation*}
The minorizations for $p=1,\infty$ now follow at once from the simple
general inequality
\begin{equation}  \label{21inf}
\|f\|_2^2 \leq \|f\|_1 \|f\|_{\infty}
\end{equation}
and the majorizations for $p=1,\infty$.

Thus we may assume $1 < p < \infty$, $p \neq 2$. The majorization follows
from the general inequality
\begin{equation}  \label{p1inf}
\|f\|_p^p \leq \|f\|_1 \|f\|_{\infty}^{p-1}
\end{equation}
and the majorizations for $1$ and $\infty$.

For the minorization, we note that if $q < 2 < r$, and if $0 < \theta < 1$
is the number with $1/2 = \theta/q + (1-\theta)/r$, then one has the general
inequality
\begin{equation}  \label{q2r}
\|f\|_2 \leq \|f\|_q^{\theta} \|f\|_{r}^{1-\theta}.
\end{equation}
If $p > 2$, the minorization follows, after a brief computation, from (\ref%
{q2r}) in the case $q = 1$, $r = p$, and the minorizations for $2$ and $1$.
If $p < 2$, the minorization follows, after a briefer computation, from (\ref%
{q2r}) in the case $q = p$, $r = \infty$, and the minorizations for $2$ and $%
\infty$. This completes the proof.

\section{Spin Besov Spaces and their Characterization}

The purpose of this section is the characterization of functional spaces by
mixed needlets. This issue was already addressed by (\cite{bkmp09}), where
the characterizations of Besov spaces by the asymptotic behaviour of spin
needlet coefficients is addressed; here we aim at an analogous goal by
focussing on mixed needlet coefficients. Most of our notations and of the
arguments to follow are classical and close to those provided by (\cite%
{bkmp09}). We start by recalling that (compare (\ref{ghio})) if $2^j \geq |s|
$,
\begin{equation*}
V_{2^{j},s}=\bigoplus_{l=|s|}^{2^{j}}\mathcal{H}_{ls}\ ,
\end{equation*}%
the space of spin functions spanned by spin spherical harmonics of degree up
to $2^{j}$. We let
\begin{equation*}
\sigma _{j}(F_s;p):=\inf_{G_s \in V_{2^{j},s}}\left\| F_s-G_s\right\|
_{L_{s}^{p}}\text{ ,}
\end{equation*}%
the error from the best approximations in that same space. The definition of
Besov spaces is then natural.

\begin{definition}
\label{spinbesdf} (\cite{bkmp09}) (Spin Besov space) We say that the spin
function $F_{s} \in L^p_s$ belongs to the Besov space of order $\left\{
p,q,r;s\right\} $ (written $F_{s}\in B_{r;s}^{pq})$ if and only if%
\begin{equation*}
\sigma _{j}(F_s;p)=\varepsilon _{j}2^{-jr},
\end{equation*}%
where $\left\{ \varepsilon _{j}\right\} \in \ell ^{q}$ and $p\geq 1$, $q,r>0,
$ $s\in \mathbf{N.}$ The associated norm is
\begin{equation*}
\left\| F_s\right\| _{B_{r;s}^{pq}}:=\left\| F_s\right\|
_{L_{s}^{p}}+\left\| \varepsilon _{j}\right\| _{\ell ^{q}}\text{ .}
\end{equation*}
\end{definition}

\textbf{Remark} Recall that $L_{j}$ is the largest integer with $%
e_{L_{j}s}\leq 2^{2(j+1)}$, and that by (\ref{ljasy}), $L_j \sim 2^j$ as $j
\to \infty$. In particular, there is a $c > 0$ such that $2^{j-c} \leq L_j
\leq 2^{j+c}$ for all $j$. Thus, if $2^{j-c} \geq |s|$,
\begin{equation*}
V_{2^{j-c},s} \subseteq V_{L_j,s} \subseteq V_{2^j,s}.
\end{equation*}
Thus if we set
\begin{equation*}
\tilde{\sigma}_{j}(F_s;p):=\inf_{G_s \in V_{L_j,s}}\left\| F_s-G_s\right\|
_{L_{s}^{p}}\text{ ,}
\end{equation*}%
it is evident that we could use the $\tilde{\sigma}_j$ in place of the $%
\sigma_j$ in Definition \ref{spinbesdf} to define the same spaces with an
equivalent norm.\newline

The following characterization is provided by (\cite{bkmp09}) and extends to
spin fiber bundles classical results on approximation spaces.

\bigskip

\begin{theorem}
\label{Besovspin} (\cite{bkmp09}) If $F_s\in L_{s}^{p},$ the following
conditions are equivalent to $F_s \in B_{r;s}^{pq}$:
\end{theorem}

\begin{enumerate}
\item
\begin{equation*}
\left\| P_{j}F_s-F_s\right\| _{L_{s}^{p}}=\varepsilon _{1j}2^{-jr}
\end{equation*}

\item
\begin{equation*}
\left\| P_{j}F_s-P_{j-1}F_s\right\| _{L_{s}^{p}}=\left\| Q_{j}F_s\right\|
_{L_{s}^{p}}=\varepsilon _{2j}2^{-jr}
\end{equation*}

\item
\begin{equation*}
\left\{ \sum_{k}\left| \beta _{jk,s}\right| ^{p}\left\| \psi _{jk,s}\right\|
_{L_{s}^{p}}^{p}\right\} ^{1/p}=\varepsilon _{3j}2^{-jr}.
\end{equation*}
\end{enumerate}

where $\left\{ \varepsilon _{1j}\right\} ,\left\{ \varepsilon _{2j}\right\}
,\left\{ \varepsilon _{3j}\right\} \in \ell ^{q},$ and $c_{p},C_{p}>0$.

\bigskip The fact that (1) implies (2) is easy, while the converse is
standard from Hardy's inequality. The fact that $F_s \in B_{r;s}^{pq}$
implies (1) follows at once from using the $\tilde{\sigma}_j$ in Definition %
\ref{spinbesdf}, while the converse follows from Corollary \ref{allgood} (a)
for the $P_j$, once one notes that for any $G_s \in V_{L_{j-1}s}$, $P_j G_s
= \phi({\textstyle\frac{\sqrt{T}}{2^{j+1}}})G_s = G_s$, so
\begin{equation*}
\|P_jF_s - F_s\|_{L_{s}^{p}} \leq \|P_j\|\|F_s - G_s\|_{L_{s}^{p}} + \|G_s -
F_s\|_{L_{s}^{p}}.
\end{equation*}
The equivalence of (2) with (3) is established by first showing that for any
$F \in L^p_s$,
\begin{equation}  \label{beseq3}
c_{p}\left\| Q_{j}F_s\right\| _{L_{s}^{p}}\leq \left\{ \sum_{k}\left| \beta
_{jk,s}\right| ^{p}\left\| \psi _{jk,s}\right\| _{L_{s}^{p}}^{p}\right\}
^{1/p}\leq C_{p}\left\{ \left\| Q_{j-1}F_s\right\| _{L_{s}^{p}}+\left\|
Q_{j}F_s\right\| _{L_{s}^{p}}+\left\| Q_{j+1}F_s\right\| _{L_{s}^{p}}\right\}
\end{equation}
We will prove a ``mixed'' analogue of (\ref{beseq3}) in Theorem \ref{Besov}
below, by a proof which is very close to the proof in \cite{bkmp09}.

\begin{remark}
In the previous Theorem, the crucial result is of course provided by (3),
which provides the characterizations of Besov classes by means of the decay
of spin needlet coefficients. This feature could be provided by many
alternative formulations; in particular, as in the mixed case, one has (\cite%
{bkmp09}) that
\begin{equation*}
\left\| \psi _{jk,s}\right\|_p \sim 2^{j(1-\frac{2}{p})}.
\end{equation*}
The previous result can hence be formulated as follows: The measurable spin
function $F_{s}$ belongs to the Besov space of order $\left\{
p,q,r;s\right\} $ if and only if
\begin{equation*}
\left[\int_{\mathbb{S}^{2}}\left| F_{s}(x)\right|^{p}dx\right]^{1/p}+%
\left[ \sum_{j}2^{qj\left\{ r+2(\frac{1}{2}-\frac{1}{p})\right\} }\left\{
\sum_{k}\left| \beta _{jk,s}\right| ^{p}\right\} ^{q/p}\right] ^{1/q}<\infty
\text{ .}
\end{equation*}
\end{remark}

Our main result in this section is to show that the mixed needlet
coefficients can play exactly the same role as the spin coefficients in the
characterization of functional spaces, despite their different mathematical
features. More precisely, we have the following alternative characterization
of Besov spaces:

\begin{theorem}
\label{Besov} The function $F_{s}\in L_{s}^{p}$ belongs to the spin Besov
space $B_{r;s}^{pq}$ if and only if
\begin{equation}
\left\{ \sum_{k}\left| \beta _{jk,s\mathcal{M}}\right| ^{p}\left\| \psi
_{jk,s\mathcal{M}}\right\| _{L_{s}^{p}}^{p}\right\} ^{1/p}=\varepsilon
_{4j}2^{-jr}.  \label{besov4}
\end{equation}%
where $\left\{ \varepsilon _{4j}\right\} \in \ell ^{q},$ and $c_{p},C_{p}>0$%
. Equivalently, $F_{s}\in B_{r;s}^{pq}$ if and only if
\begin{equation*}
\left[ \int_{\mathbb{S}^{2}}\left| F_{s}(x)\right| ^{p}dx\right] ^{1/p}+%
\left[ \sum_{j}2^{jq\left\{ r+2(\frac{1}{2}-\frac{1}{p})\right\} }\left\{
\sum_{k}\left| \beta _{jk,s\mathcal{M}}\right| ^{p}\right\} ^{q/p}\right]
^{1/q}<\infty \text{ .}
\end{equation*}
\end{theorem}

\begin{proof}
Given Theorem \ref{Besovspin} and the results established in the previous
section, the proof is rather standard and very close, for instance, to the
arguments in (\cite{bkmp09}).

By Theorem \ref{Besovspin}, it suffices to establish that for all $F \in
L^p_s$,
\begin{equation}  \label{bescrit}
c_{p}\left\| Q_{j}F_s\right\| _{L_{s}^{p}}\leq \left\{ \sum_{k}\left|
\left\langle F_{s},\psi _{jk,s\mathcal{M}}\right\rangle \right|^{p} \left\|
\psi _{jk,s\mathcal{M}}\right\| _{L_{s}^{p}}^{p}\right\} ^{1/p}\leq
C_{p}\left\{ \left\| Q_{j-1}F_s\right\| _{L_{s}^{p}}+\left\|
Q_{j}F_s\right\| _{L_{s}^{p}}+\left\| Q_{j+1}F_s\right\|
_{L_{s}^{p}}\right\} \text{.}
\end{equation}

In addition to Lemma \ref{lpnrm}, we will need the inequality
\begin{equation}
\sum_{k}\left| \psi _{jk,s\mathcal{M}}(x)\right| \leq C_{M}\sum_{k}\frac{%
C_{M}2^{j}}{\left\{ 1+2^{j}d(x,\xi _{jk})\right\} ^{M}}\leq C_{M}2^{j}\text{
.}  \label{aux2}
\end{equation}%
which follows from the properties of $\epsilon$-nets (see \cite{bkmpBer} or %
\cite{bkmp09}).\newline

To establish the rightmost inequality of (\ref{bescrit}), we note first
that, in view of the mixed analogue of (\ref{qjsix}),
\begin{equation*}
\sum_{k}\left| \left\langle F_{s},\psi _{jk,s\mathcal{M}}\right\rangle
\right| ^{p}\left\| \psi _{jk,s\mathcal{M}}\right\|
_{L_{s}^{p}}^{p}=\sum_{k}\left| \left\langle
Q_{j-1}F_{s}+Q_{j}F_{s}+Q_{j+1}F_{s},\psi _{jk,s\mathcal{M}}\right\rangle
\right| ^{p}\left\| \psi _{jk,s\mathcal{M}}\right\| _{L_{s}^{p}}^{p}\text{ .%
}
\end{equation*}%
(In fact, the sums are clearly termwise equal at least for $F_s \in
C^{\infty}_s$; the equality for general $F_s \in L^p_s$ follows by a density
argument.)\newline

The result will then follow if we can prove that, for all $G_{s}\in L_{s}^{p}
$
\begin{equation*}
\sum_{k}\left| \left\langle G_{s},\psi _{jk,s\mathcal{M}}\right\rangle
\right| ^{p}\left\| \psi _{jk,s\mathcal{M}}\right\| _{L_{s}^{p}}^{p}\leq
C^{p}\left\| G_{s}\right\| _{L_{s}^{p}}^{p}\text{ .}
\end{equation*}%
In view of Lemma \ref{lpnrm} and (\ref{aux2}), the result can be established
along exactly the same lines as in Lemma 14 of (\cite{bkmp09}); more
precisely, by Holder's inequality%
\begin{equation*}
\left| \left\langle G_{s},\psi _{jk,s\mathcal{M}}\right\rangle \right| \leq
\left\{ \int_{\mathbb{S}^{2}}\left| G_{s}(x)\right| ^{p}\left| \psi _{jk,s%
\mathcal{M}}(x)\right| dx\right\} ^{1/p}\left\{ \int_{\mathbb{S}^{2}}\left|
\psi _{jk,s\mathcal{M}}(x)\right| dx\right\} ^{1-1/p}
\end{equation*}%
whence%
\begin{eqnarray*}
\sum_{k}\left| \left\langle G_{s},\psi _{jk,s\mathcal{M}}\right\rangle
\right| ^{p}\left\| \psi _{jk,s\mathcal{M}}\right\| _{L_{s}^{p}}^{p} &\leq
&\left( \int_{\mathbb{S}^{2}}\left| G_{s}(x)\right| ^{p}\sum_{k}\left| \psi
_{jk,s\mathcal{M}}(x)\right| dx\right) \left\| \psi _{jk,s\mathcal{M}%
}\right\| _{L_{s}^{1}}^{p-1}\left\| \psi _{jk,s\mathcal{M}}\right\|
_{L_{s}^{p}}^{p} \\
&\leq &C\left\| G_{s}\right\| _{L_{s}^{p}}^{p}\text{ ,}
\end{eqnarray*}%
as desired, by Lemma \ref{lpnrm} and (\ref{aux2}). The proof of the
rightmost inequality is hence completed.

As far as the leftmost inequality of (\ref{bescrit}) is concerned, again our
arguments are very close to \cite{bkmp09}, Lemmas 15 and 16. Note first
that, by (\ref{qjtwo}) and Lemma \ref{qjm},
\begin{equation*}
Q_{j}F_{s}=\sum_{k}\left\langle F_{s},\psi _{jk,s\mathcal{M}}\right\rangle
\psi _{jk,s\mathcal{M}}\text{ ,}
\end{equation*}%
at least for $F_{s}\in L_{s}^2$. Now using Holder's inequality%
\begin{equation*}
\left\| \sum_{k}\left\langle F_{s},\psi _{jk,s\mathcal{M}}\right\rangle
\psi _{jk,s\mathcal{M}}\right\| _{L_{s}^{p}}^{p}
\end{equation*}%
\begin{eqnarray*}
&\leq &\int_{\mathbb{S}^{2}}\left\{ \sum_{k}\left| \left\langle F_{s},\psi
_{jk,s\mathcal{M}}\right\rangle \right| \left| \psi _{jk,s\mathcal{M}%
}(x)\right| ^{1/p}\left| \psi _{jk,s\mathcal{M}}(x)\right| ^{1-1/p}\right\}
^{p}dx \\
&\leq &\int_{\mathbb{S}^{2}}\sum_{k}\left| \left\langle F_{s},\psi _{jk,s%
\mathcal{M}}\right\rangle \right| ^{p}\left| \psi _{jk,s\mathcal{M}%
}(x)\right| \left\{ \sum_{k}\left| \psi _{jk,s\mathcal{M}}(x)\right|
\right\} ^{p-1}dx \\
&\leq &C2^{j(p-1)}\sum_{k}\left| \left\langle F_{s},\psi _{jk,s\mathcal{M}%
}\right\rangle \right| ^{p}\left\| \psi _{jk,s\mathcal{M}}\right\|
_{L_{s}^{1}}\leq C\sum_{k}\left| \left\langle F_{s},\psi _{jk,s\mathcal{M}%
}\right\rangle \right| ^{p}\left\| \psi _{jk,s\mathcal{M}}\right\|
_{L_{s}^{p}}^{p}\text{ ,}
\end{eqnarray*}%
again by (\ref{aux2}) and Lemma \ref{lpnrm}. This gives the leftmost
inequality of (\ref{bescrit}) for $F_{s}\in L_{s}^2$, and hence for 
$F_s \in L^p_s$ if $p \geq 2$.  If instead $p < 2$, and $%
F_{s}\in L_{s}^{p}$ is general, we take a sequence $F_{s}^{m}\in C_{s}^{\infty }$
approaching $F_{s}$ in $L_{s}^{p}$, consider that inequality with $F_{s}^{m}$
in place of $F_{s}$, and take the limsup of both sides in $m$, to obtain the
inequality for $F_{s}$, as desired.
\end{proof}

Theorem \ref{Besov} could be formulated more directly as: The section $F_{s}$
belongs to the spin Besov space $B_{r;s}^{pq}$ if and only if there exist $%
\left\{ \varepsilon _{j}\right\} \in \ell ^{q}$ such that%
\begin{equation*}
\sum_{k}\left| \beta _{jk,s\mathcal{M}}\right| ^{p}=\varepsilon
_{j}2^{-j\left\{ r+2(\frac{1}{2}-\frac{1}{p})\right\} }\text{ }.
\end{equation*}%
Combining Theorems \ref{Besov} and \ref{Besovspin}, we have the bounds%
\begin{equation*}
c\varepsilon _{2,j}2^{-jr+2j(\frac{1}{2}-\frac{1}{p})}\leq \left\{
\sum_{k}\left| \beta _{jk,s\mathcal{M}}\right| ^{p}\right\} ^{1/p},\left\{
\sum_{k}\left| \beta _{jk,s}\right| ^{p}\right\} ^{1/p}\leq C\left\{ \frac{%
\varepsilon _{2,j-1}}{2}+\varepsilon _{2,j}+2\varepsilon _{2,j+1}\right\}
2^{-jr+2j(\frac{1}{2}-\frac{1}{p})},
\end{equation*}%
where the $\ell ^{q}$ sequence $\left\{ \varepsilon _{2j}\right\} $ is such
that%
\begin{equation*}
\left\| Q_{j}F_{s}\right\| _{L_{s}^{p}}=\varepsilon _{2j}2^{-jr}.
\end{equation*}%
More explicitly, the asymptotic behaviour of the norms of needlet
coefficients is of the same order for the spin and mixed spin case, despite
the fact that the coefficients in the two cases have a rather different
nature (the $\left\{ \beta _{jk,s}\right\} $ are spin-valued, while $\left\{
\beta _{jk,s\mathcal{M}}\right\} $ are complex valued scalars). An
alternative way to formulate this conclusion is the following. Define $%
B_{r}^{pq}(\mathbb{C)}\ $\ as the Besov space of complex-valued functions on
the sphere. Then:

\begin{center}
The spin function $F_{s}$ belongs to the spin Besov space $B_{r;s}^{pq}$ if
and only if the scalar complex-valued function $f=(f_{E}+if_{M})$ belongs to
$B_{r}^{pq}(\mathbb{C)}$ , i.e.
\begin{equation*}
\left\{ F_{s}\in B_{r;s}^{pq}\right\} \Leftrightarrow \left\{
(f_{E}+if_{M})\in B_{r}^{pq}(\mathbb{C)}\right\} \text{ }.
\end{equation*}
\end{center}

Note that the complex-valued function $f$ does not correspond to the
function $g$ introduced in (\ref{nofun}), indeed for a given array of
coefficients $\left\{ a_{lm,s}=a_{lm}^{E}+ia_{lm}^{M}\right\} _{lm}$ we can
write%
\begin{equation*}
F_{s}(x)=\sum_{lm}a_{lm,s}Y_{lm,s}(x)=\sum_{lm}a_{lm,s}\sqrt{\frac{(l-s)!}{%
(l+s)!}}(\eth )^{s}Y_{lm}=(\eth )^{s}g(x)\text{ ,}
\end{equation*}%
whence%
\begin{equation*}
g(x)=\sum_{lm}\sqrt{\frac{(l-s)!}{(l+s)!}}a_{lm,s}Y_{lm}(x)\neq
f(x)=\sum_{lm}a_{lm,s}Y_{lm}(x)\text{ .}
\end{equation*}

\section{Statistical applications}

\subsection{Estimation of angular power spectra and cross-spectra}

A major asset explaining the success of needlets for the analysis of
cosmological data refers to their uncorrelation properties. More precisely,
it was shown in (\cite{bkmpAoS}), that for isotropic random fields, needlet
coefficients are asymptotically uncorrelated at any fixed angular distance
as the frequency $j$ diverges. This result was extended to the Mexican
needlet case by \cite{spalan},\cite{mayeli}, and motivated many applications
to astrophysical data, for instance (cross-)angular power spectrum
estimation (\cite{pbm06}, \cite{fay08}), detection of asymmetries (\cite%
{pietrobon1}), bispectrum estimation (\cite{ejslan}, \cite{rudjord1}, \cite%
{pietrobon2}) and many others. An analogous property was established
for spin
needlets in \cite{gelmar}; statistical techniques were then developed in %
\cite{glm}, while applications to CMB polarization data were detailed in %
\cite{ghmkp}. In this Section, we shall show how mixed needlets allow for
further applications which have great physical interest and are not feasible
by the pure spin approach, such as, for instance, estimation of
cross-spectra between scalar and spin fields.

To this aim, we shall focus on zero-mean, isotropic spin Gaussian random
fields. As discussed by \cite{gelmar},\cite{glm}, \cite{leosa} and \cite{mal}%
, the latter can be characterized by assuming that $\left\{
a_{lm}^{E},a_{lm}^{M}\right\} $ are complex-valued Gaussian random sequences
satisfying%
\begin{equation*}
Ea_{lm}^{E}=Ea_{lm}^{M}=0\text{ , }Ea_{lm}^{E}\overline{a}_{l^{\prime
}m^{\prime }}^{E}=\delta _{l}^{l^{\prime }}\delta _{m}^{m^{\prime }}C_{l}^{E}%
\text{ , }Ea_{lm}^{M}\overline{a}_{l^{\prime }m^{\prime }}^{M}=\delta
_{l}^{l^{\prime }}\delta _{m}^{m^{\prime }}C_{l}^{M},\text{ }=-l,...,l\text{
,}
\end{equation*}%
and
\begin{equation*}
a_{lm}^{E}=\overline{a}_{l,-m}^{E}\text{ ,
}a_{lm}^{M}=\overline{a}_{l,-m}^{M}\text{ .}
\end{equation*}%
For $m=0,$ $\left\{ a_{l0}^{E},a_{l0}^{M}\right\} $ are real-valued
Gaussian
with the same moments. In the cosmological literature, the sequences $%
\left\{ C_{l}^{E},C_{l}^{M}\right\} $ are known as the angular power spectra
of the $E$ and $M$ modes; clearly in the Gaussian case they encode the full
information on the dependence structure of the random field. In these area
of applications, data are collected also on a standard scalar field (the
so-called temperature of CMB radiation), which is again assumed to be
Gaussian and isotropic with angular power spectrum%
\begin{equation*}
Ea_{lm}^{T}\overline{a}_{l^{\prime }m^{\prime }}^{T}=\delta _{l}^{l^{\prime
}}\delta _{m}^{m^{\prime }}C_{l}^{T}\text{ , }Ea_{lm}^{T}\overline{a}%
_{l^{\prime }m^{\prime }}^{E}=\delta _{l}^{l^{\prime }}\delta
_{m}^{m^{\prime }}C_{l}^{TE}\text{ , }Ea_{lm}^{T}\overline{a}_{l^{\prime
}m^{\prime }}^{M}=\delta _{l}^{l^{\prime }}\delta _{m}^{m^{\prime
}}C_{l}^{TM}.
\end{equation*}%
The cross-spectra $\left\{ C_{l}^{TE},C_{l}^{TM}\right\} $ are themselves of
great physical relevance. The former is used to constrain cosmological
parameters, in particular the so-called reionization epoch, while a
detection of a non-zero value for the latter would entail an (unexpected)
violation of parity invariance at the cosmological scales. Note that, for
jointly isotropic random fields $T,E,$ we have
\begin{equation*}
C_{l}^{TE}=Ea_{lm}^{T}\overline{a}_{lm}^{E}=E\left\{ \left[ \func{Re}%
(a_{lm}^{T})+i\func{Im}(a_{lm}^{T})\right] \left[ \func{Re}(a_{lm}^{E})-i%
\func{Im}(a_{lm}^{E})\right] \right\}
\end{equation*}%
\begin{equation*}
=E\left\{ \func{Re}(a_{lm}^{T})\func{Re}(a_{lm}^{E})\right\} +E\left\{ \func{%
Im}(a_{lm}^{T})\func{Im}(a_{lm}^{E})\right\} +iE\left\{ \func{Re}(a_{lm}^{T})%
\func{Im}(a_{lm}^{E})\right\} -iE\left\{ \func{Im}(a_{lm}^{T})\func{Re}%
(a_{lm}^{E})\right\}
\end{equation*}%
\begin{equation*}
=E\left\{ \func{Re}(a_{lm}^{T})\func{Re}(a_{lm}^{E})\right\} +E\left\{ \func{%
Im}(a_{lm}^{T})\func{Im}(a_{lm}^{E})\right\} =C_{l}^{ET}\text{ ,}
\end{equation*}%
because $\left\{ \func{Re}(a_{lm}^{T}),\func{Im}(a_{lm}^{E})\right\} \overset%
{d}{=}\left\{ \func{Im}(a_{lm}^{T}),\func{Re}(a_{lm}^{E})\right\} $ by
isotropy, where $\overset{d}{=}$ denotes equality in distribution; the
latter property follows (as in \cite{bm} and \cite{mpjmva}, compare Theorem
7.2 in \cite{gelmar}) from
\begin{equation*}
\left(
\begin{array}{c}
D_{mm^{\prime }}^{l}(R)a_{lm}^{T} \\
D_{mm^{\prime }}^{l}(R)a_{lm}^{E}%
\end{array}%
\right) \overset{d}{=}\left(
\begin{array}{c}
a_{lm}^{T} \\
a_{lm}^{E}%
\end{array}%
\right) \text{ , for all }R\in SO(3)\text{ ,}
\end{equation*}%
where $\left\{ D_{m_{1}m_{2}}^{l}(R)\right\} _{m_{1},m_{2}}$ denotes as
usual the family of irreducible unitary representations of $SO(3)$ by means
of Wigner's matrices (\cite{VIK}, \cite{vmk}). It follows, in particular,
that the cross-power spectrum is always real-valued.

We shall now provide an uncorrelation result that generalizes \cite{bkmpAoS}
and \cite{gelmar}.

\begin{theorem}
Let $F_{s}$ and $T$ be jointly isotropic spin and scalar random fields
(respectively), with angular power spectra such that%
\begin{eqnarray*}
C_{l}^{E} &=&g_{E}(l)l^{-\alpha _{E}},\text{ }C_{l}^{M}=g_{M}(l)l^{-\alpha
_{M}},C_{l}^{T}=g_{T}(l)l^{-\alpha _{T}},\text{ } \\
\alpha _{E},\alpha _{M},\alpha _{T} &>&2\text{ , }\left|
g_{E}^{(i)}(u)\right| ,\left| g_{M}^{(i)}(u)\right|,\left| g_{T}^{(i)}(u)\right| \leq c_{i}u^{-i},\text{ }%
c_{i}>0\text{ , }i=0,1,2,...
\end{eqnarray*}%
Assume also that for sufficiently large $l$%
\begin{equation*}
\left| g_{E}(l)\right| ,\left| g_{M}(l)\right|,\left|
g_{T}(l)\right| >c>0\text{ .}
\end{equation*}%
Then for all $\tau >0$ there exist $C_{\tau }>0$ such that%
\begin{equation*}
\left| Corr(\beta _{jk_{1};E},\beta _{jk_{2};E})\right| ,\left|
Corr(\beta _{jk_{1};M},\beta _{jk_{2};M})\right|,\left| Corr(\beta
_{jk_{1};T},\beta _{jk_{2};T})\right| \leq \frac{C_{\tau }}{\left\{
1+2^{j}d(\xi _{jk_{1}},\xi _{jk_{2}})\right\} ^{\tau }}\text{ .}
\end{equation*}
\end{theorem}

\begin{proof}
For the coefficents of the scalar random field $T$, the proof was
provided in \cite{bkmpAoS} (see \cite{gelmar} for the extension to
the spin case). In view of the expressions provided in Section 4 for
the needlet coefficients $(\beta _{jk;E},\beta _{jk;M})$ and the
discussion following equation (\ref{ven}), the proof is identical to
the argument for the scalar case, and it is hence omitted for
brevity's sake.
\end{proof}

As we mentioned, mixed needlets allow for statistical procedures which were
unfeasible in the scalar and pure spin cases. Consider for instance the
issue of testing for a non-zero value of the cross-spectrum $C_{l}^{TE},$
which is one of the main objectives of the ESA satellite mission Planck.
Let us introduce the estimators%
\begin{equation*}
\widehat{\Gamma }_{j}^{TE}=\func{Re}\left\{ \sum_{k}\beta _{jk;E}\beta
_{jk;T}\right\} \text{ , }\widehat{\Gamma }_{j}^{TM}=\func{Re}\left\{
\sum_{k}\beta _{jk_{1};M}\beta _{jk;T}\right\} \text{ .}
\end{equation*}%
We have easily, for $A=E,M$%
\begin{eqnarray*}
E\widehat{\Gamma }_{j}^{TA} &=&E\func{Re}\left\{
\sum_{l_{1}m_{1}}\sum_{l_{2}m_{2}}b(\frac{\sqrt{e_{l_{1}s}}}{B^{j}})b(\frac{%
\sqrt{e_{l_{2}s}}}{B^{j}})a_{l_{1}m_{1};T}\overline{a}_{l_{2}m_{2};E}%
\sum_{k}Y_{l_{1}m_{1}}(\xi _{jk})\overline{Y}_{l_{2}m_{2}}(\xi _{jk})\lambda
_{jk}\right\} \\
&=&\sum_{l}b^{2}(\frac{\sqrt{e_{l_{1}s}}}{B^{j}})\frac{2l+1}{4\pi }%
C_{l}^{TA}.
\end{eqnarray*}%
Likewise
\begin{equation*}
Var\left\{ \widehat{\Gamma }_{j}^{TA}\right\} =\left\{ \sum_{l}b^{4}(\frac{%
\sqrt{e_{ls}}}{B^{j}})\frac{2l+1}{4\pi }C_{l}^{T}C_{l}^{A}\right\} +\left\{
\sum_{l}b^{4}(\frac{\sqrt{e_{ls}}}{B^{j}})\frac{2l+1}{4\pi }\left(
C_{l}^{TA}\right) ^{2}\right\} \text{ }.
\end{equation*}%
The following result is straightforward:

\begin{lemma}
\label{stoch}As $j\rightarrow \infty ,$ for $A=E,M$%
\begin{equation*}
\frac{\widehat{\Gamma }_{j}^{TA}-E\widehat{\Gamma }_{j}^{TA}}{\sqrt{%
Var\left\{ \widehat{\Gamma }_{j}^{TA}\right\} }}\rightarrow _{d}N(0,1)\text{
.}
\end{equation*}
\end{lemma}

\begin{proof}
It suffices to note that%
\begin{equation*}
\widehat{\Gamma }_{j}^{TA}=\sum_{l}b^{2}(\frac{\sqrt{e_{l_{1}s}}}{B^{j}})%
\frac{1}{4\pi }\left[ a_{l0}^{A}a_{l0}^{T}+2\sum_{m=1}^{l}\left( \func{Re}%
(a_{lm}^{A})\func{Im}(a_{lm}^{T})+\func{Re}(a_{lm}^{T})\func{Im}%
(a_{lm}^{A})\right) \right]
\end{equation*}%
and the summands satisfy all the assumptions of the classical Lindeberg-Levy
Central Limit Theorem, see \cite{gelmar} for an analogous argument.
\end{proof}

\begin{remark}
As in \cite{glm}, it is indeed possible to prove stronger results than Lemma \ref%
{stoch}, namely, it can be shown that the same limiting result holds, if the
estimator is constructed by using only the coefficients belonging to a
connected subset of $\mathbb{S}^{2}.$ This result is important for
astrophysical applications, where observations are available only on subsets
of the sphere, due to various forms of astrophysical contamination.
\end{remark}

\begin{remark}
It would be straightforward to exploit the properties of mixed needlet
coefficients for many other statistical applications. To provide an example,
it is possible to advocate estimation of the joint bispectrum of scalar and
spin random fields, along the lines of the procedures advocated for the
scalar case by \cite{ejslan} (see also \cite{m}, \cite{mptrf}, \cite{mpjmva}%
). While these extensions are straightforward from the mathematical
point of view, and hence omitted here for brevity's sake, they are
certainly of great practical importance for applications to CMB
datasets, as we discussed in the Introduction.
\end{remark}

\subsection{Spin Nonparametric Regression}

While the uncorrelation properties of needlet coefficients have
already been widely exploited in statistical inference, the
characterization provided for spin Besov spaces entails even richer
statistical opportunities which are still almost completely open for
research (see also \cite{hardle},\cite{kp2004}) for classical
results on adaptive nonparametric regression, \cite{kookim, kimkoo}
for optimal spherical deconvolution methods, \cite{bkmpAoSb} and
\cite{kpp} for some results on needlet-based shrinkage estimation
for densities on the sphere, and \cite{kim} for adaptive
nonparametric regression on vector bundles).

We envisage, in particular, applications to spin nonparametric regression by
means of mixed-needlets shrinkage, in the following sense. Consider the
regression model%
\begin{equation*}
Y_{s}(X_{i})=F_{s}(X_{i})+\varepsilon _{i}\text{ ,
}i=1,2,...,n\text{ ,}
\end{equation*}%
where $\left\{ X_{i}\right\} _{i=1,...,n}$ are (deterministic or
stochastic) locations on the sphere $X_{i}\in \mathbb{S}^{2},$
$F_{s}\in B_{r;s}^{pq}$ is a deterministic section of the spin fiber
bundle and $\left\{ \varepsilon _{i}\right\} _{i=1,...,n}$ is a
sequence of random observational errors, themselves spin $s$
variables. From the point of view of applications, we have in mind
measurements of so-called weak gravitational lensing effects (see
for instance \cite{bridle}); here, $F_{s}$ is the shear induced by
gravitational effects on the image of distant galaxies, and $\left\{
\varepsilon _{k}\right\} $ are observational errors, due for
instance to the intrinsic variability in the shape of the galaxies.
The aim is to reconstruct $F_{s}$ upon observations on $\left\{
Y_{s}(X_{i})\right\} $; this is the object of a number of ongoing
challenges, detailed for instance in \cite{bridle}. Assume $\left\{
X_{i}\right\} _{i=1,...,n}$ make up an (approximate) sequence of
cubature points; we suggest to estimate $F_{s}$ by means of the
shrinkage procedure%
\begin{equation*}
\widetilde{F}_{s}(x):=\sum_{jk}\widetilde{\beta
}_{jk,s\mathcal{M}}^{\ast }\psi _{jk,s\mathcal{M}}(x)\text{ ,}
\end{equation*}%
where%
\begin{equation*}
\widetilde{\beta }_{jk,s\mathcal{M}}^{\ast }=\widetilde{\beta }_{jk,s%
\mathcal{M}}\mathbb{I}(\left| \widetilde{\beta }_{jk,s\mathcal{M}}\right|
>ct_{n})\text{ , }t_{n}\rightarrow 0\text{ as }n\rightarrow \infty \text{
and }c>0\text{ ,}
\end{equation*}%
and
\begin{equation*}
\widetilde{\beta }_{jk,s\mathcal{M}}:=\frac{4\pi}{n}\sum_{i}Y_{s}(X_{i})\overline{\psi }%
_{jk,s\mathcal{M}}(X_{i})\simeq \int_{\mathbb{S}^{2}}F_{s}(x)%
\overline{\psi }_{jk,s\mathcal{M}}(x)dx\text{ .}
\end{equation*}%
The Besov space characterization opens the possibility of
investigating optimality properties (in the minimax sense) of the
shrinkage estimator $\widetilde{F}_{s}$
over Besov balls $B_{r;s}^{pq}(Q),$ where $Q<\infty$, defined as%
\begin{equation*}
F_{s}:\left\| F_{s}\right\| _{B_{r;s}^{pq}}<Q\text{ .}
\end{equation*}%
A full investigation of these issues will be reported elsewhere.

\end{document}